\documentclass[a4paper,12pt]{amsart}

\usepackage{graphics,graphicx}
\usepackage{tikz-cd}
\usetikzlibrary{calc}

\usepackage{xcolor}

\usepackage{enumitem}

\usepackage{booktabs}

\usepackage[a4paper]{geometry}
\usepackage[utf8]{inputenc}
\usepackage[T1]{fontenc}

\usepackage{url}

\usepackage{amssymb,amsmath,amsthm}

\usepackage[backref=page,colorlinks=true,allcolors=blue]{hyperref} %

\newtheorem{thm}{Theorem}[section]
\newtheorem{cor}[thm]{Corollary}
\newtheorem{lem}[thm]{Lemma}
\newtheorem{prop}[thm]{Proposition}
\newtheorem{lemma}[thm]{Lemma}

\newtheorem*{ack}{Acknowledgements}

\theoremstyle{definition}
\newtheorem{defn}[thm]{Definition}
\newtheorem{rem}[thm]{Remark}
\newtheorem{remark}[thm]{Remark}
\newtheorem{exa}[thm]{Example}

\DeclareMathOperator{\Ric}{Ric}

\DeclareMathOperator{\Aut}{Aut}

\DeclareMathOperator{\bary}{bar}
\DeclareMathOperator{\SGr}{SGr}

\newcommand{\red}{\mathrm{red}}
\newcommand{\scal}{\mathrm{Scal}}
\newcommand{\vol}{\mathrm{Vol}}
\newcommand{\ric}{\mathrm{Ric}}

\newcommand{\ent}{\mathbf{Ent}}
\newcommand{\Fut}{\mathbf{F}}

\newcommand{\w}{\mathrm{w}}
\renewcommand{\v}{\mathrm{v}}
\newcommand{\p}{\mathtt{p}}
\newcommand{\compatible}{\mathrm{comp}}
\newcommand{\alg}{\mathrm{alg}}

\newcommand{\T}{\mathbb{T}}
\newcommand{\bbC}{\mathbb{C}}

\newcommand{\bbZ}{\mathbb{Z}}
\newcommand{\bbP}{\mathbb{P}}
\newcommand{\bbR}{\mathbb{R}}

\newcommand{\bbS}{\mathbb{S}}
\renewcommand{\ref}{\mathrm{ref}}
\newcommand{\ka}{[\omega_0]}
\newcommand{\M}{\bf{M}}
\newcommand{\J}{\mathbf{J}}

\newcommand{\D}{\mathbf{D}}
\newcommand{\TC}{\mathbb{T}^{\mathbb{C}}}
\newcommand{\R}{\mathbb{R}}
\newcommand{\G}{\mathbb{G}}

\begin{document}

\title[Weighted numerical invariants]{Numerical invariants for weighted cscK metrics}

\author{Thibaut Delcroix}
\address{Thibaut Delcroix, Univ Montpellier, CNRS, Montpellier, France}
\email{thibaut.delcroix@umontpellier.fr}
\urladdr{https://delcroix.perso.math.cnrs.fr/}

\author{Simon Jubert}
\address{ Sorbonne Université, Université Paris Cité, CNRS, IMJ-PRG, F-75005 Paris, France}
\email{simonjubert@gmail.com}
\urladdr{https://sites.google.com/view/simon-jubert/accueil}

\date{\today}

\begin{abstract}
In K-stability, the delta invariant of a Fano variety encodes the existence of Kähler-Einstein metrics. 
We introduce a weighted analytic delta invariant, and a reduced version, that characterize the existence of weighted solitons. 
We further prove a sufficient condition of existence of weighted cscK metrics in terms of this invariant. 
We elucidate the relation between the weighted delta invariant and the greatest lower bound on the weighted Ricci curvature, called the weighted beta invariant. 
We provide a general upper bound for the weighted beta invariant in terms of moment images. 
Finally, we investigate how the geometry of semisimple principal fibrations, whose basis is not assumed to be cscK, allows to estimate their beta invariant in terms of the basis and the weighted fiber. 
Most of our statements are new even in the trivial weights settings, that is, for Kähler-Einstein and cscK metrics. 
\end{abstract}



\maketitle


\section{Introduction}

The algebraic delta invariant was introduced by Fujita and Odaka \cite{Fujita-Odaka_2018} and proved to characterize K-stability of Fano varieties in \cite{Fujita-Odaka_2018,Blum-Jonsson_2020}. 
The latter is the algebro-geometric counterpart to the existence of Kähler-Einstein metrics, thanks to the solution of the Yau-Tian-Donaldson conjecture for Fano manifolds by Chen-Donaldson-Sun \cite{CDSI, CDSII, CDSIII} and Tian \cite{Tian_2015}, and later for Fano varieties by Chi Li \cite{Li_2022}. 
The algebraic delta invariant has proved a formidable tool to effectively investigate K-stability of Fano varieties, for example via the technique to obtain lower bound developed by Abban and Zhuang \cite{Abban-Zhuang_2022}, see e.g. the book \cite{the_book}. 
Kewei Zhang introduced in \cite{Z21} an analytic variant of the algebraic delta invariant, encoding the coercivity of the Ding and Mabuchi functionals, and hence, the existence of Kähler-Einstein metrics as well. 
It was proven to coincide with the algebraic delta invariant in \cite{Zhang_2024}.  

In this article, we consider more generally the weighted (analytic) delta invariant, encoding the coercivity of the weighted entropy functional. 
The setting is as follows. 
Let \(X\) be a compact Kähler manifold. 
Let \(\Aut_{\red}(X)\) be its reduced automorphism group, that is, the connected subgroup of the group of automorphism \(\Aut(X)\) whose Lie algebra consists of vector fields which vanish somewhere on \(X\). 
Let \(\T\simeq (\bbS^1)^r\) be a compact torus acting on \(X\), factorizing through \(\Aut_{\red}(X)\). 
We consider a \(\T\)-invariant Kähler class \(\ka\) on \(X\) and a reference \(\T\)-invariant Kähler form \(\omega_0\in \ka\). 
Since the action of \(\T\) factorizes through \(\Aut_{\red}(X)\), this action is Hamiltonian, and there exists a moment map \(\mu_0\) with moment image \(\Delta=\mu_0(X)\). 
We fix a positive weight on \(\Delta\), that is, a smooth positive function \(\v:\Delta \to \bbR_{>0}\). 
For any Kähler form \(\omega\in \ka\), we consider its moment map \(\mu_{\omega}\) that satisfies \(\mu_{\omega}(X)=\Delta\), and define the \emph{weighted Ricci curvature} as
\[ \Ric_{\v}(\omega)=\Ric(\omega)-\frac{1}{2}dd^c\log(\v(\mu_\omega)) \]
This is the form associated with the Bakry-\'Emery tensor introduced early in Riemannian geometry \cite{Bakry-Emery_1985}, for a special choice of \(\T\)-invariant functions that vary consistently along a Kähler class. 

This setting allows one to consider weighted solitons, which form a large and by now well-studied generalization of Kähler-Einstein metrics (Mabuchi already considered in \cite{Mabuchi_2003} the notion of weighted Ricci curvature above and weighted solitons for a weight obtained by composition of a smooth one-variable function with a linear projection from \(\Delta\) to an integral). 
A Kähler form \(\omega\in \ka\) is a \(\v\)-soliton if
\[ \ric_{\v}(\omega)=\lambda \omega \]
for some \(\lambda\in \bbR\).
Since the Kähler class of \(\ric_{\v}(\omega)\) is \(2\pi c_1(X)\) and \(\Aut_{\red}(X)\) is trivial if \(c_1(X)\leq 0\) \cite{Carrell-Lieberman_1973}, if the weight is non-trivial we necessarily have \(\lambda>0\), \(\ka=\frac{2\pi c_1(X)}{\lambda}\) and \(X\) is Fano. 
We also consider twisted weighted solitons: given a \(\T\)-invariant Kähler form \(\theta\) on \(X\) such that \(2\pi c_1(X) = \lambda \ka +[\theta] \) for some \(\lambda\in \bbR\), we say that \(\omega\in \ka\) is a twisted weighted soliton if 
\[ \ric_{\v}(\omega)=\lambda \omega + \theta \]
Here, the underlying manifold needs not be Fano, and \(\lambda\) may be non-positive. 

Our first goal is to characterize the existence of weighted solitons and twisted weighted solitons in terms of the weighted delta invariant that we now define. 
Let \(\mathcal{K}(X,\omega_0)^{\T}\) denote the space of smooth \(\T\)-invariant strictly \(\omega_0\)-psh functions on \(X\). 
We define the weighted entropy functional on \(\mathcal{K}(X,\omega_0)^{\T}\) by 
\[ \ent_{\v}(\varphi) = \int_X \log\left(\frac{\v(\mu_{\omega_{\varphi}})\omega_{\varphi}^{[n]}}{\v(\mu_{\omega_{0}})\omega_0^{[n]}}\right) \v(\mu_{\omega_{\varphi}})\omega_{\varphi}^{[n]} \]
Consider a group  \(\G\subset \Aut_{\red}^{\T}(X)\), then we define the \emph{\(\G\)-reduced weighted delta invariant} \(\delta^{\G}_{\v}(\ka)\) as the best coercivity constant 
\[ \delta^{\G}_{\v}(\ka) = \sup \left\{\delta\in \bbR \mid \inf_{\mathcal{K}(X,\omega_0)^{\T}} \ent_{\v}-\delta (\J_{\v,\tilde{\v}}^{\omega_0})_{\G} \in \bbR  \right\} \]
where the weighted energy functional \((\J_{\v,\tilde{\v}}^{\omega_0})_{\G}\) (that we introduce in Section~\ref{sec_weighted_energy}) is a functional that measures coercivity modulo additive constants and action of the group \(\G\). It is a weighted analogue of the functional \(I-J\) considered usually. 
When \(\G\) is trivial, we omit the group in the notation and call \(\delta_{\v}(\ka)\) the weighted delta invariant of \(\ka\) (in this case, it was introduced by Rubinstein, Tian and Zhang in \cite{RTZ}). 

\begin{thm}
\label{thm_intro_soliton}
Assume that \(X\) is Fano, that \(\v\) is log concave \footnote{In \cite{HL}, which is used in our proofs, there is a gap in the proof of smoothness of weak weighted solitons for non log concave weights. J. Han kindly informed us that he has a fix for this gap, but we add the assumption of log concavity for now since this argument is not widely available.} and that the image of \(\T\) in \(\Aut_{\red}(X)\) is a maximal compact torus, then \(X\) admits a \(\v\)-soliton if and only if the weighted Futaki invariant of \(X\) vanishes and 
\(\delta_{\v}^{\TC}(2\pi c_1(X)) > 1 \), where \(\TC\) denotes the complexified torus. 
\end{thm}

It is essential to consider a \emph{reduced} version of the delta invariant here since there is not a unique \(\T\)-invariant \(\v\)-soliton, but a full \(\Aut_{\red}^{\T}(X)\)-orbit of \(\v\)-solitons, so the weighted Mabuchi functional cannot be coercive on the space of Kähler metrics. 
The weighted Futaki invariant used in the statement above is an obstruction to the existence of canonical metrics first introduced by Futaki for Kähler-Einstein metrics \cite{Futaki_1983}. In fact, the weighted Mabuchi functional is affine along \(\Aut_{\red}^{\T}(X)\) orbits in \(\mathcal{K}(X,\omega_0)^{\T}\), and the weighted Futaki invariant encodes its slope.  
To the best of our knowledge, no such result has been considered before in the literature in terms of an analytic invariant, even in the Kähler-Einstein setting. 
From the algebraic point of view, the characterization in terms of an algebraic reduced delta invariant (in the unweighted setting) of reduced uniform K-stability, in turn equivalent to K-polystability and existence of Kähler-Einstein metrics is due to \cite{Xu-Zhuang_2020,Li_2022,LXZ_2022}.

We further introduce the greatest weighted Ricci lower bound 
\[ \beta_\v(\ka) := \sup \{\beta \in \bbR\mid \exists 0 < \omega \in \ka, \T\text{-invariant st}, \ric_\v(\omega)-\beta\omega>0 \}. \]
The unweighted greatest Ricci lower bound for Fano manifolds was first considered in early work of Tian \cite{Tian_1992} and Rubinstein \cite{Rubinstein_2009}, but it gained much more interest following the work of Székelyhidi \cite{Szekelyhidi_2011}, who showed that it coincides with the maximal time of existence of a solution in the continuity method approach to Kähler-Einstein metrics. 
This point of view allowed Zhang and Berman, Boucksom, and Jonsson to relate it with the analytic or algebraic delta invariant \cite[Proposition~3.8]{Z21}, \cite[Theorem~C]{BBJ21}, and earlier, by Cheltsov, Rubinstein and Zhang \cite[Appendix]{CRZ_2019} in the Fano case. 
We prove the weighted version of this statement. 
For this, consider \(s(\ka):=\sup\{ s \mid 2\pi c_1(X)-s\ka >0 \}\).

\begin{thm}
\label{thm_intro_beta_delta}
For a log concave weight \(\v\), the weighted beta invariant is characterized as 
\begin{align*}
 \beta_{\v} (\ka) & = \min\{s(\ka),\delta_{\v}(\ka)\}  \\ 
  &= \sup \{t<s(\ka)\mid  \forall 0<\theta\in (2\pi c_1(X)-t\ka), \T-\text{invariant}, \\
& \phantom{= \sup \{t<s(\ka) \quad} \exists \varphi \in \mathcal{K}(X,\omega_0)^{\T} \text{ st } \ric_{\v}(\omega_{\varphi})=t\omega_{\varphi}+\theta \} 
\end{align*}
\end{thm}

This result was stated without a detailed proof in \cite[Proposition~6.15]{RTZ}. 
The greatest weighted Ricci lower bound had already been considered earlier in the literature, under the name of the greatest Bakry-Emery-Ricci lower bound, and computed for certain Fano manifolds with many symmetries (e.g. \cite{DGSW_2018}, \cite{Delcroix-Hultgren_2021}). 
In particular, for these varieties the previous statement allows to compute a lower bound for \(\delta_{\v}(\ka)\). 

By elementary considerations on the moment maps, we show that the formula for the greatest weighted Ricci lower bound of a toric manifold provides a general upper bound for the greatest weighted Ricci lower bound outside of the toric setting. 

\begin{thm}
\label{thm_intro_general_upper}
Let \(\bary_{\v}(\ka) = \frac{\int_X \mu_{\omega} \v(\mu_{\omega})\omega^{[n]}}{\int_X \v(\mu_{\omega})\omega^{[n]}} \in \mathfrak{t}^*\) be the barycenter of \(\Delta\) with respect to the weighted Duistermaat-Heckman measure, then 
\[ \beta_{\v}(\ka) \leq \sup \{\beta < s(\ka) \mid -\beta\bary_{\v}(\ka) \in \Delta_{2\pi c_1(X)-\beta\ka} \} \]
where the moment image \(\Delta_{2\pi c_1(X)-\beta\ka}\) is chosen according to the choice of moment image \(\Delta\) of \(\ka\), and the canonical moment image of \(2\pi c_1(X)\). 
\end{thm}

To our knowledge, this statement is new even for \(\v \equiv 1\). 
An alternative way to obtain an upper bound, valid for \(\delta_\v\), would be to consider the asymptotics of the functionals along the \(\Aut_{\red}^{\T}(X)\)-orbits. However, note that neither the entropy nor the weighted energy functional are linear along these orbits, so their slopes are not so easy to express. 
Székelyhidi used this approach in the Fano case in \cite[Proposition~7]{Szekelyhidi_2011}, in which case non-discrete automorphism group imposes \(1\geq \delta(2\pi c_1(X))=\beta(2\pi c_1(X))\), and the bound obtained seems to be equivalent to the one we obtain here (at least in the applications known to the authors).

Let \(\w:\Delta\to \bbR\) be another smooth function, not assumed to take positive values. We now consider the weighted Mabuchi functional \(\mathbf{M}_{\v,\w}\) (see Section~\ref{sec_weighted_Mab} for the definition) whose minimizers are the weighted cscK metrics introduced by Lahdili \cite{Lah19}. 
In the unweighted case, a strategy developed by Song and Weinkove \cite{Wei03,Weinkove_2006,SW08} inspired by the introduction of the \(\J\)-flow by Donaldson \cite{Donaldson_1999} and early work of Chen \cite{Chen_2000,Chen_2004}, 
allows to obtain a sufficient condition of coercivity of the Mabuchi functional \emph{via} the study of the \(\J\)-flow. 
Several variants have been studied throughout the years e.g. \cite{LSY,Li-Shi_2016,Dervan_2016,Zheng_2015,JSS_2019,Sjostrom_Dyrefelt_2022,Z21,To_2023}. 
Following this general strategy, we prove a sufficient condition for the existence of weighted cscK metrics in terms of the reduced weighted delta invariant. 

\begin{thm}
\label{thm_intro_exi:delta}
Assume that the weighted Futaki invariant $\Fut_{\v,\w}$ vanishes. 
Assume that for some \(\varepsilon>0\), \(\theta_{\varepsilon}:=\delta_{\varepsilon}\ka-2\pi c_1(X)\) is Kähler, where \(\delta_{\varepsilon}:=\delta_{\v}^{\TC}(\ka)-\varepsilon\). 
Consider the function \(\check{\w}_{\varepsilon} : \Delta\times \Delta_{\theta_{\varepsilon}} \to \bbR\) defined by 
\[ \check{\w}_{\varepsilon}(x,y) = \delta_{\varepsilon}n-\frac{1}{2}\w(x)+\langle d\log \v(x),\delta_{\varepsilon}x-y\rangle. \]
Assume that either: 
\begin{enumerate}
    \item the weight \(\v\) satisfies \(1+\langle d\log(\v)(x),x-x'\rangle \geq 0\) for all \(x\), \(x'\in \Delta\), that \(\inf \check{\w}_{\varepsilon} +(n-1)(s(\ka)-\delta_{\varepsilon})>0\) and that \(\forall y\in \Delta_{\theta_{\varepsilon}}\), \(x\mapsto \check{\w}_{\varepsilon}(x,y)\) is convex, or
    \item \(X=\bbP^1\), \(\ka=2\pi c_1(\bbP^1)\) and \(\inf \check{\w}_{\varepsilon}>0\). 
\end{enumerate}
Then the weighted Mabuchi functional $\M_{\v,\w}$ is $\TC$-coercive.
If \(\v\) is furthermore log concave, then there exists a weighted cscK metric. 
\end{thm}

Note again that this statement is new in the case when \(\v\equiv 1\) since we work with the \emph{reduced} delta invariant again. 
To pass from coercivity of the weighted Mabuchi functional to existence of weighted cscK metrics, we rely on the results of Di Nezza, Lahdili and the first author \cite{DJL, DJLb} and Han and Liu \cite{HanLiu} independently generalizing in the weighted setting the breakthrough of Chen and Cheng on the variational approach to cscK metrics \cite{Chen-Cheng_2021a,Chen-Cheng_2021b}.  
The strategy of the proof of coercivity relies first on comparing the weighted Mabuchi functional with a well-chosen weighted energy functional, using the weighted delta invariant. 
The other ingredient is to obtain a lower bound on this weighted energy functional via the study of weighted $\mathbf{J}$-equation, which requires to adapt the elliptic approach of Li--Shi--Yao~\cite{LSY} and Li--Shi \cite{Li-Shi_2016} to the general weighted setting.  

We stress that our criterion above should be seen only as a first step in this direction, and we hope that it will motivate the study of the weighted $\mathbf{J}$-equation
\begin{equation}
\Lambda_{\varphi,\v}(\theta) = \frac{1}{2}\w(\mu_{\varphi})
\end{equation}
and of the weighted \(\J\)-flow 
\begin{equation}
\frac{\partial \varphi}{\partial t} = \frac{1}{2}\w(\mu_{\varphi})-\Lambda_{\varphi,\v}(\theta) 
\end{equation}
where we refer to Section~\ref{sec_weighted_trace} for the definition of the weighted trace \(\Lambda_{\varphi,\v}(\theta)\). 
For example, our sufficient criterion of existence of a solution to the $\mathbf{J}$-equation depends very much on the Kähler form involved, rather than its Kähler class. 
It is natural to hope for a result along the lines of that of Chen \cite{Chen_2021} in the unweighted case, at least in the toric setting with the approach of \cite{Collins-Szekelyhidi_2017,Takahashi}. 
Deriving best coercivity constants for weighted energy functionals in the spirit of the work of Sjöström~Dyrefelt \cite{Sjostrom_Dyrefelt_2020} is also a natural direction to explore. 

Note that the first assumption on \(\v\) in the first situation of Theorem~\ref{thm_intro_exi:delta} says that \(\v\) has small variations with respect to its absolute values. 
We expect that this assumption could be removed by following a parabolic approach rather than an elliptic approach. 
The extreme case is when \(\v\) is constant, in which case the conditions are much simpler to check. In this very particular case, the reader may be interested in the recent works \cite{Hisamoto-Nakamura,Apostolov-Lahdili-Nitta} for sufficient conditions of existence of weighted cscK metrics with a completely different approach. 

The convexity assumption in the first situation of Theorem~\ref{thm_intro_exi:delta} is quite strong, but satisfied for example if \(\log \v\) is affine and \(\w\) is convex. In particular, this setting includes the weights for cscK metrics, extremal Kähler metrics, Kähler-Ricci solitons and combinations. Our criterion is furthermore much more powerful than the one obtained in \cite{Li-Shi_2016} since we use here the reduced delta invariant which is the best possible constant to use in this approach.  
It would be interesting to study the weighted \(\J\)-equation from the point of view of subsolutions to possibly remove the assumption, but it does not seem to fit directly in the range of equations considered in \cite{Szekelyhidi_2018}, so a detailed study would be needed. 

Finally, let us comment on the second situation of Theorem~\ref{thm_intro_exi:delta}. 
Since the manifold \(X\) is \(\bbP^1\), the reader may think that the statement is trivial, but far from it, this situation applies to the case of Calabi ansatz, or more generally simple principal fibrations, which correspond to the fibrations introduced below, with one-dimensional fiber. 

As we have highlighted above, most of our statements are new even in the unweighted setting. 
Furthermore, for non-trivial weights, we obtain applications to the unweighted setting by considering semisimple principal fibrations. 
A semisimple principal fibration \(Y\) is a fiber bundle with structure group a torus, with fiber a compact complex manifold \(X\), and some additional assumptions (see Section~\ref{sec_sspf}). There is a fiberwise action of a torus \(\T\) on \(Y\), and in favorable cases, the geometry of \(Y\) reduces to the weighted geometry of \(X\) for a weight associated with the construction. 
These were first studied in \cite{ACGT_2011} when the fiber is toric, and in \cite{AJL} for general fiber. 
In contrast to these previous works, we \emph{do not require} the basis to be cscK. 
One of our goal is to elucidate more precisely when the geometry of \(Y\) reduces to that of the weighted fiber. 

We study this \emph{via} the \emph{compatible weighted beta invariant} \(\beta^{\compatible}_{\v}(Y,[\tilde{\omega}_0])\), defined as the supremum of all real numbers \(\beta\) such that there exists a Kähler form \(\tilde{\omega}\in [\tilde{\omega}_0]\) which is compatible (with the semisimple principal fibration structure, see Definition~\ref{defn_comp}) and such that \(\ric_{\v}(\tilde{\omega})>\beta \tilde{\omega}\). 
There is an obvious inequality 
\(\beta^{\compatible}_{\v}(Y,[\tilde{\omega}_0]) \leq \beta_{\v}(Y,[\tilde{\omega}_0])\) that we show to be strict by examples (see Example~\ref{exa_not_fiber}). 
We compute \(\beta^{\compatible}_{\v}(Y,[\tilde{\omega}_0])\) in terms of the basis and the weighted fiber in general. 
Instead of reproducing the general statement (Theorem~\ref{thm:beta_comp}) in this introduction, let us focus on a special case. 

We say that a semisimple principal fibration \(Y\) is \emph{compatibly Fano} if there exists a compatible Kähler metric \(\tilde{\omega}\in 2\pi c_1(Y)\) with \(\ric(\tilde{\omega})\) Kähler. 
For such fibrations, we prove the following expression for their weighted invariant. We refer to Section~\ref{sec_fibrations} for details on the notations used, but let us note that for the two terms in the minimum, one depends only on the weighted fiber, while the other depends only on the basis factors. 

\begin{thm}
\label{thm_intro_comp_Fano}
If \(Y\) is a compatibly Fano, then 
\[ \beta_{\v}^{\compatible}(Y,2\pi c_1(Y)) = \min \left\{ \beta_{\v\p}(X,2\pi c_1(X)), \min_a \frac{c_a \beta(B_a,2\pi c_1(B_a))+\inf_{\Delta_{2\pi c_1(X)}}p_a}{c_a+\inf_{\Delta_{2\pi c_1(X)}}p_a} \right\} \] 
If furthermore the fiber is toric, and the minimum above is achieved by \(\beta_{\v\p}(X,2\pi c_1(X))\), then 
\[ \beta_{\v}(Y,2\pi c_1(Y)) = \beta_{\v}^{\compatible}(Y,2\pi c_1(Y)) = \beta_{\v\p}(X,2\pi c_1(X)) \]
\end{thm}

The latter condition shows that the weighted fiber encodes the geometry of the fibration well beyond the case when the basis is Kähler-Einstein. 
It further provides a quantitative bound on where this relation starts to fail. 
In Example~\ref{exa_not_fiber} we provide an example where the minimum is not achieved by \(\beta_{\v\p}(X,2\pi c_1(X))\). 
In fact, this specific example belongs to a much more restricted class of fibrations, the projective Fano \(\bbP^1\)-bundles considered by Zhang and Zhou in \cite{ZZ}. 
We provide a detailed comparison of our results with theirs in Section~\ref{sec_comp_ZZ}. 
It is important to note that, contrary to their result, we do not fully compute the beta invariant but only the compatible beta invariant which is a lower bound. 
On the other hand, our results provides a geometric interpretation of a lower bound which is often sharp, it does detects when the beta invariant is computed by the weighted fiber, and it provides generalization both to the weighted setting, and to much more general fibration constructions. 

The notion of compatibly Fano itself is a further evidence of how the weighted fiber may not encode well the geometry of the fibration. 
We characterize this property in terms of the beta invariant of the basis and show by example that it is stronger than being Fano in Example~\ref{exa_not_comp_Fano}.
To find such examples, we consider \(\bbP^1\)-bundles over the odd symplectic Grassmannians \cite{Mihai_2007}, and exploit the latter family of varieties' property to have high index, and low beta invariant as the dimension grows \cite{Hwang-Kim-Park_2023}.

Let us finally mention a simple application of our result, which allows us to derive properties of the basis (and the fiber) from properties of the total space, unlike the usual results on semisimple principal fibration. 

\begin{cor} 
\label{cor_intro_Fano}
If \(Y\) is a compatibly Fano semisimple principal fibration and satisfies \(\beta^{\compatible}(2\pi c_1(Y))\geq 1\), then \(Y\) is K-semistable, the basis \(B\) is K-semistable and the fiber \(X\) is weighted K-semistable. 
\end{cor}

The article is organized as follows. 
In Section~\ref{sec_recall_weighted}, we recall the setting of weighted Kähler geometry, as well as some of the basic results and definitions needed. 
We introduce the weighted delta invariant in Section~\ref{sec_delta_and_soliton}, where we also prove the criterion for existence of weighted solitons: Theorem~\ref{thm_intro_soliton} is Theorem~\ref{t:exi:sol:delta}. 
We introduce the weighted beta invariant and prove Theorem~\ref{thm_intro_beta_delta} (=Theorem~\ref{thm_beta_delta}) in Section~\ref{sec_beta_delta}. 
Section~\ref{sec_J_cscK} is devoted to the proof of Theorem~\ref{thm_intro_exi:delta} (=Theorem~\ref{t:exi:delta}), and the proof of a priori estimates for solutions of the weighted $\mathbf{J}$-equation occupies a significant portion of the article. 
The final Section~\ref{sec_fibrations} deals with semisimple principal fibrations and proves Theorem~\ref{thm_intro_comp_Fano} and Corollary~\ref{cor_intro_Fano} and provides the examples alluded to in this introduction. 

\begin{ack}
The second author is grateful to Abdellah Lahdili for numerous enlightening conversations. He would also thank Chi Li, Jiyuan Han and Chenxi Yin for useful discussions. Both authors thank Tat Dat Tô for useful comments and exchanges on this work. The second author is funded by the ERC SiGMA - 101125012 (PI: Eleonora Di Nezza). 
The first author is partially funded by ANR-21-CE40-0011 JCJC project MARGE.
\end{ack}

\section{Weighted Kähler geometry}
\label{sec_recall_weighted}

\subsection{Weighted setting}
Let \( X \) be a compact Kähler manifold of dimension \(n\), and let \(\Aut(X)\) denote its automorphism group. 
The reduced automorphism group \( \Aut_{\mathrm{red}}(X) \) is the subgroup of \(\Aut(X)\) whose Lie algebra consists of vector fields which vanish somewhere on \(X\). 
Let \(\T\simeq (\bbS^1)^r\) be a (non-trivial) compact torus acting holomorphically on \(X\), and such that its action factors through \( \Aut_{\mathrm{red}}(X) \).
We denote by \(\mathfrak{t}\) the Lie algebra of \(\T\), and by \(\mathfrak{t}^*\) its dual, both isomorphic to \(\bbR^r\). 

As a consequence, for any \(\T\)-invariant closed real \((1,1)\)-form \(\theta\) on \(X\), one can associate a  \emph{moment map} \(\mu:X\to \mathfrak{t}^*\) (in a generalized sense) satisfying 
\[ d\langle \mu, \xi \rangle = -\iota_{\xi}\theta \]
for all \(\xi\in \mathfrak{t}\).
When \(\theta\) is a Kähler form, hence also a symplectic form, this is a moment map in the usual sense. 
Moment maps behave linearly under operations on the forms: if \(\mu_1\) and \(\mu_2\) are moment maps for \(\theta_1\) and \(\theta_2\), and \(t\in \bbR\), then \(t\mu_1+\mu_2\) is a moment map for \(t\theta_1+\theta_2\). 

The moment map \(\mu\) is only defined up to translation by an element of \(\mathfrak{t}^*\), but there are some natural choices in certain cases. 
First, if \(\varphi: X\to \bbR\) is a smooth function on \(X\), then \(d^c\varphi\) is a moment map for \(dd^c\varphi\), where \(d^c\varphi\) is interpreted as a map to \(\mathfrak{t}^*\), such that \(\langle d^c\varphi ,\xi\rangle = d^c\varphi(\xi) \) where we identify an element \(\xi\in \mathfrak{t}\) with the vector field it generates on \(X\). 
Second, if \(\theta\in 2\pi c_1(X)\) (whether or not \(X\) is Fano, since we consider \((1,1)\)-forms that may not be Kähler), then the natural lift of the action of \(\T\) to the anticanonical line bundle \(K_X^{-1}\) induces a natural choice of moment map. 
Whenever a special case as above occurs, we will use the naturally normalized moment map. 

We fix a \( \T \)-invariant Kähler form \( \omega_0 \) on \( X \) with moment map \( \mu_0 : X \to \mathfrak{t}^* \). We let \( \Delta := \mu_0(X) \)  be the moment image of \( (X,\omega_0,\T) \)  which is a convex polytope \cite{Atiyah_1982,Guillemin-Sternberg_1982}.  We consider the space of \( \mathbb{T} \)-invariant \( \omega_0 \)-relative Kähler potentials  
\[ 
\mathcal{K}(X,\omega_0)^\mathbb{T} := \left\{ \varphi \in \mathcal{C}^{\infty}(X)^\mathbb{T} \mid \omega_\varphi := \omega_0 + dd^c \varphi > 0 \right\}.
\]
For any Kähler potential \( \varphi \in \mathcal{K}(X,\omega_0)^\mathbb{T} \), we denote by 
\[
\mu_\varphi = \mu_0 + d^c\varphi,
\]
the moment map associated with \(\omega_{\varphi}\), in such a way that its image \(\mu_{\varphi}(X)=\Delta\) does not depend on \( \varphi \), but only on the class \([\omega_0]\).   

In view of the previous observations, given this fixed choice of moment map \(\mu\) for \(\omega_0\), there is a natural normalization of the moment map for any form in a linear combination \(t_1[\omega_0]+t_2 2\pi c_1(X)\) with \(t_1,t_2\in \bbR\), and we denote its moment image by 
\(\Delta_{t_1[\omega_0]+t_2 2\pi c_1(X)}\) if the class is Kähler. 

We fix a function \( \v \in \mathcal{C}^{\infty}(\Delta, \mathbb{R}_{>0}) \), called a \textit{weight function}. 

\subsection{Weighted curvatures and weighted trace}
\label{sec_weighted_trace}

The central idea of weighted Kähler geometry is to replace the volume form \(\omega_{\varphi}^n\) with the volume form \(\v(\mu_{\varphi})\omega_{\varphi}^n\). 
Note that the pushforward of this volume form by the moment map \(\mu_{\varphi}\) is independent of \(\varphi\), this is the measure with potential \(\v\) with respect to the Duistermaat-Heckman measure. 
We denote by 
\[ 
\vol_{\v} := \int_{\Delta} \v DH = \int_X \v(\mu_{\varphi})\omega_{\varphi}^n 
\] 
the total mass of this measure. 

The \textit{weighted Ricci form} is defined by 
\begin{equation}{\label{v:ricci}}
    \Ric_\v(\omega_\varphi) := \Ric(\omega_\varphi) -\frac{1}{2}dd^c \log \v(\mu_\varphi).
\end{equation}
For the Ricci form, we use the convention of \cite{Gau}. In particular, the cohomology class of the Ricci form of any K\"ahler metric $\omega$ is given by $[\Ric(\omega)]= 2\pi c_1(X)$, where $c_1(X)$ is the first Chern class of the anticanonical line bundle $K_X^{-1}$.

Following \cite{Lah19, DJL}, for any $(1,1)$-form $\rho$ \emph{with moment map} $\mu_{\rho}$, we define the $\v$-\textit{weighted trace} of $\rho$ with respect to \(\omega_{\varphi}\) by
\begin{equation*}
    \Lambda_{\varphi,\v}(\rho):= \Lambda_{\varphi} (\rho) + \langle d\log(\v)(\mu_{\varphi}), \mu_{\rho} \rangle 
\end{equation*}
\noindent for any $\varphi \in \mathcal{K}(X,\omega_0)^\T$,
where \(\Lambda_{\varphi} (\rho):=n\frac{\rho\wedge\omega_\varphi^{n-1}}{\omega_\varphi^{n}}\) denotes the usual trace of \(\rho\) with respect to \(\omega_{\varphi}\). 

According to \cite{Lah19, DJL}, for any K\"ahler metric we introduce the $\v$-\textit{weighted laplacian} of $\omega_\varphi$ acting on $\T$-invariant smooth function on $X$

\begin{equation*}
    \Delta_{\varphi,\v}(f):= \Lambda_{\varphi,\v}dd^cf,
\end{equation*}
A standard computation shows that

\begin{equation*}
    \Delta_{\varphi,\v}(f)=d^{\star}_{\varphi,\v}df,
\end{equation*}
where $d^{\star}_{\varphi,\v}$ is the formal adjoint operator of $d$ with respect to the \textit{weighted volume form} $\v(\mu_{\varphi})\omega_{\varphi}^{[n]}$, see \cite[Appendix A]{DJL} for more details. 
We will make use of the resulting weighted variant of integration by part, which reads 
\begin{equation}
\label{eq_wibp}
\int_X f \Delta_{\varphi,\v}(g) \v(\mu_{\varphi})\omega_{\varphi}^{n} = \int_X g \Delta_{\varphi,\v}(f) \v(\mu_{\varphi})\omega_{\varphi}^{n}
\end{equation}

Finally, the \emph{$\v$-weighted scalar curvature} of a $\T$-invariant K\"ahler metric $\omega_\varphi \in [\omega_0]$ is defined by
\begin{equation}{\label{def:scalv2}}
 {\scal_\v(\omega_\varphi)}:= 2 \Lambda_{\varphi,\v}(\Ric_\v(\omega_\varphi)),
\end{equation}

\begin{remark}
 In \cite{Lah19, AJL, DJL}, the weighted scalar curvature is defined by ${\scal^{\mathrm{Lah}}_\v(\omega_\varphi)}:= 2 \v(\mu_\varphi)\Lambda_{\varphi,\v}(\Ric_\v(\omega_\varphi))$. Hence, our definition of the weighted scalar curvature is related to the one use in \cite{Lah19, AJL, DJL} by  ${\scal^{\mathrm{Lah}}_\v(\omega_\varphi)}= \v(\mu_\varphi) {\scal_\v(\omega_\varphi)}$.
\end{remark} 

\subsection{Weighted constant scalar curvature Kähler metrics}

Given a second weight function $\w \in {C}^{\infty}(\Delta,\R)$ (note that we do not impose positive values for \(\w\), in contrast to \(\v\)), a $\T$-invariant K\"ahler metric $\omega_\varphi$ is called \emph{$(\v,\w)$-weighted cscK} if its $\v$-weighted scalar curvature satisfies
\begin{equation}
\label{eqn_wcscK}
    \scal_\v(\omega_\varphi)=\w(\mu_\varphi).
\end{equation}
We mention that, from our definition of the weighted scalar curvature \eqref{def:scalv2}, the choice of $\w$ for encoding certain geometric situations varies by a factor of $\v$ comparing to \cite{Lah19, AJL}.

When \(X\) is Fano and \(\omega_0\in 2\pi c_1(X)\), a particular case of interest in this paper is that of a \(\v\)-\textit{soliton}, which is, by definition, a K\"ahler metric \(\omega_\varphi\) such that
\begin{equation}{\label{v-sol-def}}
    \Ric_\v(\omega_\varphi)=\omega_\varphi.
\end{equation}
As proven in \cite[Lemma~2.2]{AJL}, \(\omega_{\varphi}\) is a \(\v\)-soliton if and only if it is a \((\v,\tilde{\v})\)-weighted cscK metric, with 
\begin{equation}
\label{defn_tildev}
\tilde{\v}(x) = 2(n+\langle d\log \v (x),x\rangle)
\end{equation}

\subsection{Weighted Mabuchi functional and the weighted Futaki invariant}
\label{sec_weighted_Mab}

As a variational approach to weighted cscK metrics, the weighted Mabuchi energy was introduced in \cite{Lah19} and is defined by its variation 
\begin{equation}
    d_{\varphi}\left(\mathbf{M}_{\v,\w}\right)(\dot{\varphi})=-\int_X \dot{\varphi}\big(\scal_\v(\omega_\varphi)-\w(\mu_\varphi)\big)\frac{\v(\mu_\varphi)\omega_\varphi^{n}}{\vol_{\v}}, \quad \mathbf{M}_{\v,\w}(0)=0.
\end{equation} 
If \(\omega_{\varphi}\) is a weighted cscK metric, then 
\begin{equation}{\label{norma:weight}}
    \int_X \Big(\scal_\v(\omega_{\varphi}) - \w(\mu_{\varphi})\Big) \v(\mu_\varphi)\omega_{\varphi}^{n}=0,
\end{equation}
But the left-hand side is independent of \(\varphi\in \mathcal{K}(X,\omega_0)^\mathbb{T}\), so provides a necessary condition for existence of cscK metrics. 

Going further, we define the weighted Futaki invariant $\mathbf{F}_{\v,\w} : \mathrm{Aff}(\Delta) \longrightarrow \mathbb{R}$ by
\begin{equation}{\label{d:weit:ft}}
\mathbf{F}_{\v,\w}(\ell):=\int_X \Big(\scal_\v(\omega) - \w(\mu_{\omega})\Big) \ell(\mu_{\omega})\v(\mu_\omega)\omega^{[n]}.
\end{equation}

The following statement is well known in the unweighted case (it was one of the motivations of Mabuchi to introduce the K-energy \cite{Mabuchi_1986}) and follows by considering the restriction of the weighted Mabuchi functional to orbits of the complexified torus \(\TC\) (see \cite[Lemma 4.6]{DJLb} for details). 

\begin{lemma}{\label{fut:van}}
Assume that \(\v\) and \(\w\) satisfy condition~\ref{norma:weight}. 
Then the weighted Mabuchi energy $\M_{\v,\w}$ is $\TC$-invariant if and only if the weighted Futaki invariant $\mathbf{F}_{\v,\w}$ vanish.
More precisely, if \(\xi\in \mathfrak{t}\) and \((\varphi_t)_{t\in \bbR}\) is a family in \(\mathcal{K}(X,\omega_0)^\mathbb{T}\) such that for all \(t\), \(\omega_{\varphi_t}=(e^{t\xi})^*\omega_{\varphi}\), then 
\begin{equation}
\label{linear_Mab}
\mathbf{M}_{\v,\w}(\varphi_t)= t \Fut_{\v,\w}(\xi) + \mathbf{M}_{\v,\w}(\varphi) 
\end{equation}
\end{lemma}

\section{Weighted delta invariant of a K\"ahler class} 
\label{sec_delta_and_soliton}

\subsection{Weighted energy functionals and coercivity}
\label{sec_weighted_energy}

Fix weights \(\v\), \(\w\) as in the previous section, and a closed real \(\T\)-invariant \((1,1)\)-form \(\rho\) with moment map \(\mu_{\rho}\). 

\begin{defn}
\label{def:j:v}
The weighted (twisted) energy functional \(\J_{\v,\w}^{\rho} : \mathcal{K}(X,\omega_0)^\T \to \bbR\) is defined by 
\begin{equation}{\label{def:j}}
\J_{\v,\w}^{\rho}(\varphi) = \int_0^1\int_X  \varphi \big(  2\Lambda_{s\varphi,\v}(\rho) - \w(\mu_{s\varphi})\big) \frac{\v(\mu_{s\varphi})\omega_{s\varphi}^{n}}{\vol_{\v}}   \mathop{ds}.
\end{equation}
\end{defn}

One useful property of this family of functionals is linearity in \((\w,\rho)\): 
\begin{equation}
\label{J-linearity}
\J_{v,t\w_1+\w_2}^{t\rho_1+\rho_2} = t\J_{v,\w_1}^{\rho_1} + \J_{v,\w_2}^{\rho_2}
\end{equation}
Furthermore, when \(\rho_1\) and \(\rho_2\) are in the same cohomology class, the difference \(\J_{\v,\w}^{\rho_1}-\J_{\v,\w}^{\rho_2}\) is bounded. 
More precisely, we have the following. 

\begin{lemma}
\label{lem_J_ddc}
Let \(f:X\to \bbR\) be a smooth function. 
Then 
\[ \J_{\v,0}^{dd^c f}(\varphi) = 2 \left(\int_X f \frac{\v(\mu_{\varphi})\omega_{\varphi}^{n}}{\vol_{\v}} - \int_X f \frac{\v(\mu_0)\omega_0^{n}}{\vol_{\v}}\right) \]
In particular, it is bounded: 
\[ \lvert \J_{\v,0}^{dd^c f}(\varphi) \rvert \leq 4 \lVert f \rVert_{C^0} \]
\end{lemma}

\begin{proof}
By definition, we have 
\begin{align*}
\J_{\v,0}^{dd^c f}(\varphi) & 
= \int_0^1\int_X  2 \varphi \Lambda_{s\varphi,\v}(dd^c f) \frac{\v(\mu_{s\varphi})\omega_{s\varphi}^{n}}{\vol_{\v}}   \mathop{ds} \\ 
& = \int_0^1\int_X  2 \varphi \Delta_{s\varphi,\v}(f) \frac{\v(\mu_{s\varphi})\omega_{s\varphi}^{n}}{\vol_{\v}}   \mathop{ds} \\ 
& = \int_0^1\int_X  2 f \Delta_{s\varphi,\v}(\varphi) \frac{\v(\mu_{s\varphi})\omega_{s\varphi}^{n}}{\vol_{\v}}   \mathop{ds} \\ 
\end{align*}
by weighted integration by parts~\eqref{eq_wibp}. 
Now a standard computation (see \cite[Appendix~A]{DJL}) shows that 
\[ \int_0^1 \Delta_{s\varphi,\v}(\varphi) \frac{\v(\mu_{s\varphi})\omega_{s\varphi}^{n}}{\vol_{\v}}   \mathop{ds} = \frac{\v(\mu_{\varphi})\omega_{\varphi}^{n}-\v(\mu_0)\omega_0^{n}}{\vol_{\v}} \]
The statements follows. 
\end{proof}

We say that a functional \(\mathbf{F} : \mathcal{K}(X,\omega_0)^\T \longrightarrow \mathbb{R}\) is \emph{invariant under addition of a constant} if for any \(\varphi\in \mathcal{K}(X,\omega_0)^\T\) and any \(c\in \bbR\), \(\mathbf{F}(\varphi+c)=\mathbf{F}(\varphi)\). 
The functional \(\J_{\v,\w}^{\rho}\) is invariant under addition of a constant if and only if 
\[ 
\int_X 2\Lambda_{\varphi,\v}(\rho)\v(\mu_{\varphi})\omega_{\varphi}^n = 
\int_X \w(\mu_{\varphi})\v(\mu_{\varphi})\omega_{\varphi}^n
\]
Note that both sides of the condition above are independent of \(\varphi\). 
This follows directly from pushforward by \(\mu_{\varphi}\) for the right-hand side, and this is proven in \cite[Lemma 2]{Lah19} for the left-hand side. 
For example, the functional \(\J_{\v,\tilde{\v}}^{\omega_0}\) where \(\tilde{\v}\) is defined by~\eqref{defn_tildev} is invariant under addition of a constant by \cite[Lemma~3.1]{DJLb}. 

We use this particular case \(\J_{\v,\tilde{\v}}^{\omega_0}\) of the functional \(\J_{\v,\w}^{\rho}\) to encode coercivity on the space of Kähler potentials, and coercivity up to a group action. 
Let \(\G\) be a complex Lie group in the connected
component of the identity of the centralizer \(\Aut_\red^\T(X)\) of \(\T\) in \(\Aut_{\red}(X)\). 
We will mainly be interested in the cases \(\G=\{1\}\) and \(\G=\T^{\bbC}\). 

\begin{defn}
 We say that a functional \(\mathbf{F} : \mathcal{K}(X,\omega_0)^\T \longrightarrow \mathbb{R}\) is $\G$-coercive  if there exists $C$, $C'$ positive constants such that for any $\varphi\in \mathcal{K}(X,\omega_0)^\T$
\begin{equation*}
    \mathbf{F}(\varphi) \geq C (\J^{\omega_0}_{\v,\tilde{\v}})_{\G}(\varphi)-C',
\end{equation*}
where $(\J^{\omega_0}_{\v,\tilde{\v}})_{\G}(\varphi):=\inf_{\gamma \in \G}\J^{\omega_0}_{\v,\tilde{\v}}(\varphi_\gamma)$ where \(\varphi_{\gamma}\) denotes an element of \(\mathcal{K}(X,\omega_0)^\T\) such that \(\gamma^*\omega_{\varphi}=\omega_{\varphi_{\gamma}}\). When $\G:=\{ 1 \}$, we say that $\mathbf{F}$ is  coercive.
\end{defn}
Note that, if \(\mathbf{F}\) is \(\G\)-coercive, then it is invariant under addition. 

This measure of coercivity will be the correct choice if one wants to recover weighted versions of the properties known for the unweighted delta invariant. 
But we first check that it is consistent with the more often used measure of coercivity via the \(J\)-functional. 
We recall that it is defined by 
\begin{equation}
\label{defn_J_standard}
J(\varphi) = \int_0^1 \int_X \varphi \frac{\omega_0^n-\omega_{s\varphi}^n}{\vol_{1}} \mathop{ds} 
\end{equation}
and correspondingly \(J_{\G}(\varphi)= \inf_{\gamma \in \G}J(\varphi_\gamma)\). 

\begin{prop}
\label{prop_coercivity}
A functional  \(\mathbf{F} : \mathcal{K}(X,\omega_0)^\T \longrightarrow \mathbb{R}\) is $\G$-coercive  if and only if there exists $C$, $C'$ positive constants such that for any $\varphi\in \mathcal{K}(X,\omega_0)^\T$
\begin{equation*}
    \mathbf{F}(\varphi) \geq C J_{\G}(\varphi)-C'
\end{equation*}
\end{prop}

\begin{proof}
This follows directly from the comparison 
\[  \frac{1}{n}\frac{\inf_{\Delta}(\v)}{\sup_{\Delta}(\v)} J(\varphi) \leq \J_{\v,\tilde{\v}}^{\omega_0}(\varphi) \leq n \frac{\sup_{\Delta}(\v)}{\inf_{\Delta}(\v)} J(\varphi) \]
which in turn follows from  \cite[Lemma~4.1 and proof of Lemma~4.2]{DJLb}, \cite[Lemma~6.4]{AJL} (see also \cite[(2.37), (2.38)]{HL}).
\end{proof}

\subsection{Weighted entropy and weighted reduced delta invariant}

We define the \emph{weighted entropy} \(\ent_{\v}(\varphi)\) of a \(\T\)-invariant Kähler form \(\omega_{\varphi}\) in \([\omega_0]\) as the entropy of the probability measure defined by \(\v(\mu_{\varphi})\omega_{\varphi}^n\) with respect to that defined by \(\v(\mu_0)\omega_0^n\): 
\begin{equation}{\label{w:ent}}
    \ent_\v(\varphi):= \int_X \log\left( \frac{\v(\mu_\varphi)\omega_\varphi^n}{\v(\mu_0)\omega_0^n}\right)\frac{\v(\mu_\varphi)\omega_\varphi^{n}}{\vol_\v} .
\end{equation}

\begin{defn}
For \(\G\subset \Aut_{\red}^{\T}(X)\), the \emph{$\v$-weighted \(\G\)-reduced delta invariant} of the Kähler class \(\ka\) is  
\begin{equation}{\label{d:delata:v:a}}
    \delta^{\G}_{\v}(\ka) := \sup \left\{\delta\in \bbR \mid \inf_{\mathcal{K}(X,\omega_0)^{\T}} \left(\ent_{\v}-\delta (\J_{\v,\tilde{\v}}^{\omega_0})_{\G}\right) \in \bbR  \right\}
\end{equation}
\end{defn}

When $\G=\{1\}$, we simply write $\delta^{\{1\}}_\v=\delta_\v$ and we refer to it as \textit{weighted delta invariant}. 
Applying Jensen's inequality and comparing $\delta^\G_{\v}$ with Tian's alpha invariant \cite{Tia87} shows that 
\begin{equation*}
     \delta^\G_{\v}(\ka)>0,
\end{equation*}
see \cite[Proposition~6.2]{Z21} for the detailed proof in the unweighted case. 
Also, we can easily see from \eqref{d:delata:v:a} that we have the following scaling property for \(t>0\).
 \begin{equation}{\label{scal:prop}}
\delta^{\mathbb{G}}_\v(\ka)=t \delta^{\mathbb{G}}_{\v(\cdot /t)}(t\ka) 
\end{equation}

As proved by Lahdili in \cite[Theorem~5]{Lah19}, the weighted Mabuchi functional admits a \emph{Chen-Tian decomposition} as a difference of a weighted entropy term and a weighted energy term:
\begin{equation}
\label{Chen-Tian}
\mathbf{M}_{\v,\w}(\varphi)=\mathbf{Ent}_\v(\varphi)  -\mathbf{J}^{\mathrm{Ric}(\omega_0)}_{\v,\w}(\varphi)
\end{equation}
We will use this decomposition repeatedly in the course of the paper to relate the delta invariant with existence of canonical Kähler metrics.  

\subsection{Weighted reduced delta invariant and existence of weighted solitons}

In this section, we will prove the following criterion of existence of weighted solitons in terms of the reduced delta invariant. 

\begin{thm}
\label{t:exi:sol:delta}
Assume that \(X\) is Fano, that \(\v\) is log concave and that the image of \(\T\) in \(\Aut_{\red}(X)\) is a maximal compact torus, then \(X\) admits a \(\v\)-soliton if and only if the weighted Futaki invariant of \(X\) vanishes and 
\(\delta_{\v}^{\TC}(2\pi c_1(X)) > 1 \), where \(\TC\) denotes the complexified torus. 
\end{thm}

The first, essential step in the proof is provided by Han and Li's characterization of the existence of \(\v\)-soliton in terms of the weighted Mabuchi functional. 

\begin{thm}{\cite[Theorem~3.6]{HL}}
\label{t:exi:sol:w}
Assume that \(X\) is Fano, that \(\v\) is log concave and that \(\T\subset \Aut_{\red}(X)\) is a maximal compact torus. Then, there exists a $\v$-soliton in $2\pi c_1(X)$ if and only if the weighted Mabuchi energy $\M_{\v,\tilde{\v}}$ is $\TC$-coercive.
\end{thm}
 
\begin{proof}[Proof of Theorem~\ref{t:exi:sol:delta}]
Consider the Ricci potential $h_{\omega_0}$ of $\omega_0$ defined by 
\begin{equation}{\label{ricci:pot}}
    \ric(\omega_0)=\omega_0+\frac{1}{2}dd^c h_{\omega_0}, \quad \quad \quad \int_X h_{\omega_0} \v(\mu_0)\omega_0^{n}=0.
\end{equation}
We claim that 
\begin{equation}
    \mathbf{M}_{\v,\tilde{\v}}(\varphi)= \ent_\v(\varphi) - \J_{\v,\tilde{\v}}^{\omega_0}(\varphi) - \frac{1}{2} \int_X h_{\omega_0}\frac{\v(\mu_\varphi)\omega_\varphi^{n}}{\vol_{\v}}.
\end{equation}
Indeed, by the Chen-Tian formula~\eqref{Chen-Tian}, then by linearity~\eqref{J-linearity} we have 
\begin{align*} 
\M_{\v,\tilde{\v}}(\varphi) & 
= \ent_\v(\varphi)  -\J^{\mathrm{Ric}(\omega_0)}_{\v,\tilde{\v}}(\varphi) \\& 
= \ent_\v(\varphi)  -\J^{\omega_0}_{\v,\tilde{\v}}(\varphi) - \frac{1}{2}\J_{\v,0}^{dd^c h_{\omega_0}}(\varphi) 
\end{align*}
Using Lemma~\ref{lem_J_ddc}, and the normalization of \(h_{\omega_0}\) yields the claim. 

Let us now prove the theorem. 
Assume first that \(X\) admits a \(\v\)-soliton. 
Then by Theorem~\ref{t:exi:sol:w}, the functional \(\M_{\v,\tilde{\v}}\) is \(\TC\)-coercive, that is, by Proposition~\ref{prop_coercivity}, there exists \(\varepsilon >0\) and \(C\in \bbR\) such that for all \(\varphi\in \mathcal{K}(X,\omega_0)^\T\), 
\[ \M_{\v,\tilde{\v}}(\varphi) \geq \varepsilon (\J_{\v,\tilde{\v}}^{\omega_0})_{\TC}(\varphi) -C \]
By our claim, we have 
\begin{align*} 
\ent_{\v}(\varphi) & 
\geq \mathbf{M}_{\v,\tilde{\v}}(\varphi) + \mathbf{J}_{\v,\tilde{\v}}^{\omega_0}(\varphi) -\frac{\lVert h_{\theta}\rVert_{C^0}}{2} 
\geq (1+\varepsilon)(\mathbf{J}_{\v,\tilde{\v}}^{\omega_0})_{\TC}(\varphi) -C-\frac{\lVert h_{\theta}\rVert_{C^0}}{2}
\end{align*}
hence \(\delta_{\v}^{\TC}(2\pi c_1(X)) \geq 1+\varepsilon >1\). 
Furthermore, since \(M_{\v,\w}\) is \(\TC\)-coercive, it is in particular bounded from below, hence \(\Fut_{\v,\tilde{\v}}\equiv 0\) by Lemma~\ref{fut:van}.   

Conversely, assume that \(\Fut_{\v,\tilde{\v}}\equiv 0\) and \(\delta_{\v}^{\TC}(2\pi c_1(X)) > 1 \). 
By  Theorem~\ref{t:exi:sol:w}, it suffices to show that \(\M_{\v,\tilde{\v}}\) is \(\TC\)-coercive. 
Let \(\varphi\in \mathcal{K}(X,\omega_0)^\T\). 
By Proposition~\ref{prop_coercivity}, the restriction of the functional \(\J_{\v,\tilde{\v}}^{\omega_0}\) to the \(\TC\)-orbit of \(\varphi\) is a proper function (up to addition of a constant to the entry) on a finite dimensional vector space, hence it admits a minimum (see \cite[Lemma 6.2]{BD87} for details in the classical case, and \cite[lemma 11]{Lah20} for the weighted case). 
So we may and do choose $\gamma_0 \in \TC$ such that \(\gamma_0^*\omega_{\varphi}=\omega_{\varphi_{\gamma_0}}\) and 
\begin{equation*}
     (\J_{\v,\tilde{\v}}^{\omega_0})_{\TC}(\varphi)=  (\J_{\v,\tilde{\v}}^{\omega_0})( \varphi_{\gamma_0}),
\end{equation*}
By Lemma~\ref{fut:van}, \(\M_{\v,\tilde{\v}}\) is \(\TC\)-invariant. 
We can thus write 
\begin{align*}
\M_{\v,\tilde{\v}}(\varphi)  
& = \M_{\v,\tilde{\v}}( \varphi_{\gamma_0}) \\
& \geq \ent_{\v}(\varphi_{\gamma_0}) - \J_{\v,\tilde{\v}}^{\omega_0}(\varphi_{\gamma_0}) - \frac{1}{2}\lVert h_{\theta}\rVert_{C^0} \\
& \geq (\delta_{\v}^{\TC}(2\pi c_1(X))-\varepsilon-1)(\J_{\v,\tilde{\v}}^{\omega_0})_{\TC}(\varphi) - \frac{1}{2}\lVert h_{\theta}\rVert_{C^0}
\end{align*}
for any \(\varepsilon > 0\). 
Since \(\delta_{\v}^{\TC}(2\pi c_1(X)) > 1 \), choosing \(\varepsilon\) small enough shows that \(\M_{\v,\tilde{\v}}\) is \(\TC\)-coercive. 
\end{proof}

\section{Greatest lower bound on weighted Ricci curvature}{\label{s:soliton}}
\label{sec_beta_delta}

\subsection{Relation between weighted delta invariant and weighted beta invariant}

We consider the greatest lower bound on the weighted Ricci curvature for \(\T\)-invariant Kähler metrics in the fixed cohomology class \([\omega_0]\). 

\begin{defn}
The weighted beta invariant of \([\omega_0]\) is 
\[ \beta_\v(\ka) := \sup \{\beta \in \bbR\mid \exists 0 < \omega \in \ka, \T\text{-invariant st}, \ric_\v(\omega)-\beta\omega>0 \}. \]
\end{defn}

It follows from the definition that \(\beta_{\v}\) has a scaling property analogous to that of the weighted delta invariant: for \(t>0\), 
\begin{equation} 
\label{eqn:beta_multiple}
\beta_\v(\ka)=t \beta_{\v(\cdot /t)}(t\ka).
\end{equation}
 We also define the K\"ahler threshold 
\begin{equation*}
    s(\ka):=\sup\{ s \in \mathbb{R} \text{ | } 2\pi c_1(X)-s\ka >0 \}.
\end{equation*}
which satisfies the scaling property \(s(\ka)=t s(t\ka)\). 

Our first goal in this section is to prove the following result. 

\begin{thm}
\label{thm_beta_delta}
\label{prop:beta_vs_delta}
For a log concave weight \(\v\), the weighted beta invariant is characterized as 
\begin{align*}
 \beta_{\v} (\ka) & = \min\{s(\ka),\delta_{\v}(\ka)\}  \\ 
  &= \sup \{t<s(\ka)\mid  \forall 0<\theta\in (2\pi c_1(X)-t\ka), \T\text{-invariant}, \\
& \phantom{= \sup \{t<s(\ka) \quad} \exists \varphi \in \mathcal{K}(X,\omega_0)^{\T} \text{ st } \ric_{\v}(\omega_{\varphi})=t\omega_{\varphi}+\theta \} 
\end{align*}
\end{thm}

As direct corollary (in fact, this is used in the proof), we can characterize the existence of twisted weighted solitons in terms of the delta invariant, when the twisting form is Kähler. 

\begin{cor}
\label{cor_twisted_solitons}
\label{c:exi:sol}
Assume that \(\theta\in 2\pi c_1(X)-\ka\) is a \(\T\)-invariant Kähler form, and that \(\v\) is log concave. Then there exists \(\varphi\in \mathcal{K}(X,\omega_0)^{\T}\) such that 
\begin{equation}
\label{eqn_twisted_soliton}
    \ric_{\v}(\omega_{\varphi})=\omega_{\varphi}+\theta
\end{equation}
if and only if \(\delta_{\v}(\ka)>1\). 
\end{cor}

This result is known in the twisted K\"ahler--Einstein case, i.e. when $\v=1$, see \cite[Theorem 1.5]{Z21}. In \cite[Theorem 6.13]{RTZ}, an analogous result is shown in finite dimensional Bergman space in term of the \textit{quantized weighted delta invariant}.

Another corollary is that, when \(s(\ka)\leq 0\), then \(\beta_{\v}(\ka) = s(\ka)\) since \(\delta_{\v}(\ka)\) is always positive. 

\subsection{Proof of Theorem~\ref{thm_beta_delta}}

We roughly follow the argument presented in \cite[Section 7.3]{BBJ21} in the unweighted setting. 
It is obvious by definition of \(\beta_{\v}(\ka)\) and \(s(\ka)\), that \(\beta_{\v}(\ka)\leq s(\ka)\). 
Let \(t<s([\omega_0])\), so that \(2\pi c_1(X)-t\ka\) is a Kähler class. 
We have to show that \(t<\delta_{\v}(\ka)\) if and only if, for any \(\theta\in  2\pi c_1(X)-t\ka\) a \(\T\)-invariant Kähler form, there exists a solution to the twisted weighted Kähler-Ricci soliton equation 
\begin{equation}
\label{eqn_twisted}
\ric_{\v}(\omega_{\varphi})=t\omega_{\varphi}+\theta 
\end{equation}
We treat the cases differently according to the sign of \(t\). 

\subsubsection{Case when \(t>0\)}

We assume that \(0<t<\delta_{\v}(\ka)\). 
Note that 
\[ \ric_{\v(\cdot/t)}(t\omega) = \ric_{\v}(\omega) \]
so that, equation~\eqref{eqn_twisted} is equivalent to 
\[ \ric_{\v(\cdot/t)}(t\omega_{\varphi}) = t\omega_{\varphi}+\theta \]
that is, to the existence of a solution to \(\ric_{\v(\cdot/t)}(\hat{\omega})=\hat{\omega}+\theta\) in the Kähler class \([\hat{\omega}_0]=t[\omega_0]\). 
In view of the scaling property~\eqref{scal:prop} of the weighted delta invariant,
it suffices now to prove Corollary~\ref{cor_twisted_solitons}, so to simplify the notations, we omit the \(t\) from now on. 

Building on \cite{Berman-Witt_Nystrom_2014, BBGZ, BBEGZ19}, Han and Li prove in \cite[Theorem~3.6]{HL} that the existence of a solution to~\eqref{eqn_twisted_soliton} is equivalent to coercivity (in our sense thanks to Proposition~\ref{prop_coercivity}) of the \emph{twisted weighted Mabuchi functional} \(\mathbf{M}_{\v,\theta}\) on \(\mathcal{K}(X,[\omega_0])^\T\). 
Up to an irrelevant constant, this functional (defined in \cite[Definition~2.10]{HL}) may be written in our notations as 
\begin{equation}{\label{twi:mab}}
    \mathbf{M}_{\v,\theta}(\varphi)= \ent_{\v}(\varphi) - \J^{\omega_0}_{\v,\tilde{\v}}(\varphi) - \int_X h_{\theta}\frac{\v(\mu_\varphi)\omega_\varphi^{n}}{\vol_{\v}},
\end{equation}
where \(h_{\theta}\) is the twisted Ricci potential of \(\omega_0\) with respect to \(\theta\), defined by 
\begin{equation}
    \ric(\omega_0) - \omega_0 - \theta =dd^c h_{\theta} \qquad \qquad \int_X h_{\theta} \v(\mu_0)\omega_0^{[n]}=0.
\end{equation}
In fact, from the same arguments as in the proof of Theorem~\ref{t:exi:sol:delta}, we have \(\mathbf{M}_{\v,\theta} = \mathbf{M}_{\v,\tilde{\v}} - \mathbf{J}_{\v,0}^{\theta}\). 

It follows directly from the definition of \(\delta_{\v}(\ka)\) that for any \(\varepsilon >0\), 
\[ \mathbf{M}_{\v,\theta}(\varphi) \geq (\delta_{\v}(\ka)-\varepsilon-1)\J^{\omega_0}_{\v,\tilde{\v}}(\varphi) - \lVert h_{\theta} \rVert_{C^0} \] 
hence it is coercive if \(\delta_{\v}(\ka)>1\). 
Conversely, if \(\mathbf{M}_{\v,\theta}\) is coercive, then there exists \(\varepsilon >0\) such that 
\[ \mathbf{M}_{\v,\theta}(\varphi) \geq \varepsilon\J^{\omega_0}_{\v,\tilde{\v}}(\varphi) \]
hence 
\[ \ent_{\v}(\varphi) 
\geq \mathbf{M}_{\v,\theta}(\varphi) + \J^{\omega_0}_{\v,\tilde{\v}}(\varphi) - \lVert h_{\theta} \rVert_{C^0} 
\geq (\varepsilon+1)\J^{\omega_0}_{\v,\tilde{\v}}(\varphi) -  \lVert h_{\theta} \rVert_{C^0} \]
hence \(\delta_{\v}(\ka)\geq 1+\varepsilon >1\). 

\subsubsection{Case when \(t=0\)}
\label{sec_weighted_CY}

In view of the previous remarks, we want to show that there always exists a solution to the equation \(\ric_{\v}(\omega_{\varphi})=\theta\).
Solving the equation \(\ric_{\v}(\omega_{\varphi})=\theta\) amounts to a weighted version of Calabi-Yau theorem. 
The proof is essentially found in the literature: existence of a weak solution follows from \cite[Theorem~1.2]{Berman-Witt_Nystrom_2014}, while regularity follows from the arguments in \cite[Section~4]{HL}. 

\subsubsection{Case when \(t<0\)}

We now want to solve the equation 
\(\ric_{\v}(\omega_{\varphi})=t\omega_{\varphi}+\theta\)
with \(t<0\). 
As it is the case with the Kähler-Einstein equation, this case is the easiest. 
We did not find a reference in the literature, so we quickly explain a proof. 
By scaling, we may again restrict to the case of the equation 
\(\ric_{\v}(\omega_{\varphi})=-\omega_{\varphi}+\theta\).

Consider again the twisted Ricci potential \(h_{\theta}\), defined now by 
\begin{equation}
    \ric(\omega_0) + \omega_0 - \theta =dd^c h_{\theta} \qquad \qquad \int_X h_{\theta} \v(\mu_0)\omega_0^{[n]}=0.
\end{equation}
Consider the twisted weighted Ding functional defined by 
\begin{equation}{\label{ding-fonc}}
        \D^{-}_{\v,\theta} (\varphi) := -\J^{0}_{\v,-1}(\varphi) + \frac{1}{2}\log\left( \int_X e^{2 \varphi+h_\theta} \frac{\omega_0^{n}}{\vol_{1}} \right),
\end{equation}
Let us prove that it is coercive, following arguments in 
\cite[Theorem 11.8 (i)]{GZ}. 

Since \(\D^{-}_{\v,\theta}\) is invariant under addition of a constant, up to replacing \(\varphi\) with \(\varphi-\sup_X\varphi\), we may as well assume that the potential \(\varphi\) satisfies \(\sup_X\varphi=0\). 
By compactness of the set of such normalized \(\omega_0\)-psh functions, there exists a constant \(C\) such that 
\[ \frac{1}{2}\log\left( \int_X e^{2 \varphi+h_\theta} \frac{\omega_0^{n}}{\vol_{1}} \right) \geq -C \]
For the remaining term, we have 
\begin{align*}
-\J^{0}_{\v,-1}(\varphi) & = \int_0^1\int_X (-\varphi)\frac{\v(\mu_{s\varphi})\omega_{s\varphi}^n}{\vol_{\v}} \mathop{ds} \\
& \geq \frac{\inf_{\Delta} \v \vol_{1}}{\vol_{\v}} \int_0^1\int_X (-\varphi)\frac{\omega_{s\varphi}^n}{\vol_{1}} \mathop{ds}\\
& \geq \frac{\inf_{\Delta} \v \vol_{1}}{\vol_{\v}} \int_0^1\int_X (-\varphi)\frac{\omega_{s\varphi}^n-\omega_0^n}{\vol_{1}} \mathop{ds}\\
& \geq \frac{\inf_{\Delta} \v \vol_{1}}{\vol_{\v}} J(\varphi)
\end{align*}
by the equation~\eqref{defn_J_standard} defining \(J\) (we remind that we assumed \(\sup_X\varphi=0\) hence \(-\varphi\geq 0\) throughout). 
By invariance under addition of a constant of both the functionals \(\D^{-}_{\v,\theta}\) and \(J\), we have shown
\[ \D^{-}_{\v,\theta} (\varphi) \geq \frac{\inf_{\Delta} \v \vol_{1}}{\vol_{\v}} J(\varphi) -C \]
that is, the functional \(\D^{-}_{\v,\theta}\) is coercive. 
The end of the proof, from coercivity to the existence of a smooth solution to the weighted Monge-Ampère equation is not different from the case treated in \cite{HL}. 
\qed

\subsection{A general upper bound from the moment images}

For \(\beta<s(\ka)\), the class \(2\pi c_1(X)-\beta\ka\) is Kähler. Let \(\Delta_{2\pi c_1(X)-\beta\ka}\) denote the moment image of this Kähler class, which is a convex polytope. 
Let also \(\bary_{\v}(\ka)\) denote the barycenter of \(\Delta\) with respect to the weight \(\v\) and the Duistermaat-Heckman measure, alternatively, 
\[ \bary_{\v}(\ka) = \frac{\int_X \mu_{\omega} \v(\mu_{\omega})\omega^{n}}{\vol_{\v}} \in \mathfrak{t}^*\]

\begin{thm}
\label{thm_general_upper}
We have 
\[ \beta_{\v}(\ka) \leq \sup \{\beta < s(\ka) \mid -\beta\bary_{\v}(\ka) \in \Delta_{2\pi c_1(X)-\beta\ka} \} \]
\end{thm}

\begin{proof}
Assume that \(\beta < \beta_{\v}(\ka)\). Then there exists Kähler forms \(\omega\in \ka\) and \(\theta\in 2\pi c_1(X)-\beta \ka\) such that 
\( \ric_{\v}(\omega)-\beta\omega = \theta \). 
In particular, at the level of moment maps we have 
\[ \mu_{\ric_{\v}(\omega)}-\beta\mu_{\omega} = \mu_{\theta}. \]
We integrate this equality with respect to the probability measure \(\frac{\v(\mu_{\omega})\omega^{n}}{\vol_{\v}}\) over \(X\) to get 
\[ \int_X \mu_{\ric_{\v}(\omega)}\frac{\v(\mu_{\omega})\omega^{n}}{\vol_{\v}} - \beta \bary_{\v}(\ka) = 
\int_X \mu_{\theta}\frac{\v(\mu_{\omega})\omega^{n}}{\vol_{\v}} \]
By standard computations (see e.g. \cite[Lemma~5]{Lah19}), we have 
\[ \mu_{\ric_\v(\omega)} = -\frac{1}{2}\Delta_{\omega,\v}\mu_{\omega}. \] 
Hence by the weighted integration by parts formula~\eqref{eq_wibp}, we have  
\[ \int_X \mu_{\ric_{\v}(\omega)}\v(\mu_{\omega})\omega^{n} = 0 \]
Furthermore, since we integrate with respect to a probability measure, 
\[ \int_X \mu_{\theta}\frac{\v(\mu_{\omega})\omega^{n}}{\vol_{\v}} \in \Delta_{[\theta]}=\Delta_{2\pi c_1(X)-\beta\ka}.\]
This finishes the proof. 
\end{proof}

To get a more tractable upper bound, we consider the case when \(X\) is Fano and \(\ka=2\pi c_1(X)\). 
Then \(2\pi c_1(X)-\beta\ka = (1-\beta) 2\pi c_1(X)\), so that \(\Delta_{2\pi c_1(X)-\beta\ka} = (1-\beta)\Delta_{2\pi c_1(X)}\). 

\begin{cor}Suppose that $X$ is Fano. Then
\[ \beta_{\v}(2\pi c_1(X)) \leq \sup \left\{\beta < 1 \mid \frac{-\beta}{1-\beta} \bary_{\v}(2\pi c_1(X)) \in \Delta_{2\pi c_1(X)} \right\}. \]
\end{cor}

This upper bound is reminiscent of the formula for the greatest Ricci lower bound for toric manifolds: it is in fact a somewhat direct generalization of the proof for the upper bound on toric manifolds that one can find e.g. in \cite{Li_2011}. 
In another direction, this upper bound was generalized to group compactifications in \cite{Delcroix_2017}, and in later works for other varieties with a lot of symmetries \cite{Yao_2017,Golota_2020}. 
It would be interesting to find a similar general upper bound for arbitrary groups actions.

\section{Existence of weighted cscK metrics via the weighted $\mathbf{J}$-equation}
\label{sec_J_cscK}

\subsection{Existence of weighted cscK metrics}

Let 
$\hat{\v} : \Delta \times \Delta \longrightarrow \R$  be defined by
\begin{equation}
\label{hatv}
\hat{\v}(x,x')= \langle d\log(\v)(x),x-x'\rangle 
\end{equation}
Let \(\varepsilon>0\) and set \(\delta_{\varepsilon}=\delta_{\v}^{\TC}(\ka)-\varepsilon\), \(\theta_{\varepsilon} = \delta_{\varepsilon}\ka-2\pi c_1(X)\) to simplify notations. 
If \(\theta_{\varepsilon}\) is Kähler, we define a function \(\check{\w}_{\varepsilon}:\Delta\times \Delta_{\theta_{\varepsilon}}\to \bbR\) by  
\begin{equation}
\label{checkw}
\check{\w}_{\varepsilon}(x,y) = \delta_{\varepsilon}n-\frac{1}{2}\w(x)+\langle d\log \v(x),\delta_{\varepsilon}x-y\rangle
\end{equation}

Our goal in this section is to prove the following sufficient condition of coercivity of the Mabuchi functional. 

\begin{thm}{\label{t:exi:delta}}
Assume that 
\begin{enumerate}[label=(\roman*)]
    \item the weighted Futaki invariant $\Fut_{\v,\w}$ vanishes,
    \item \(\theta_{\varepsilon}\) is Kähler,
    \item \(\inf \check{\w}_{\varepsilon} +(n-1)(s(\ka)-\delta_{\varepsilon})>0\),
    \item \(1+ \inf(\hat{\v})>0\), and 
    \item \(\forall y\in \Delta_{\theta_{\varepsilon}}\), \(x\mapsto \check{\w}_{\varepsilon}(x,y)\) is convex. 
\end{enumerate}
Then the weighted Mabuchi functional $\M_{\v,\w}$ is $\TC$-coercive.
\end{thm}

In the case of dimension one, which is far from trivial and applies to higher dimensional geometric situations (see Section~\ref{sec_fibrations}), we can do much better. 

\begin{thm}
\label{thm_cscK_P1}
Assume that \(X=\bbP^1\),  \(\T=\bbS^1\) and  \(\ka=2\pi c_1(X)\). 
Assume that the weighted Futaki invariant $\Fut_{\v,\w}$ vanishes, that \(\delta_{\varepsilon}>1\) and \(\inf \check{\w}_{\varepsilon}>0\). 
Then the weighted Mabuchi functional $\M_{\v,\w}$ is $\TC$-coercive.
\end{thm}

As explained in the introduction, we will follow the strategy initiated by Song and Weinkove \cite{SW08}, in the variant proposed by Li--Shi--Yao \cite[Theorem 1.1]{LSY} in the unweighted case. 
Thanks to the characterization of existence of weighted cscK metrics in terms of the coercivity of the weighted Mabuchi functional proved by Di Nezza, the second author and Lahdili in \cite{DJL,DJLb}, and independently by Han and Liu in \cite{HL}, we have as a corollary a sufficient condition of existence of weighted cscK metrics when \(\v\) is log concave. 
Note that the converse direction (from existence to coercivity) is shown for $(\v,\w)$-extremal K\"ahler metric in \cite{AJL}.

\begin{cor}
Under the same hypothesis as in Theorem~\ref{t:exi:delta} or in Theorem~\ref{thm_cscK_P1}, if \(\v\) is furthermore log concave, then there exists a $(\v,\w)$-cscK metric in $\ka$.
\end{cor}

The first key argument to obtain the theorem is to use the Chen-Tian formula~\eqref{Chen-Tian} to reduce to a lower bound on a weighted energy functional. 

\begin{lem}
\label{lem_reduce_to_J}
If the weighted Futaki invariant \(\Fut_{\v,\w}\) vanishes, and \(\J_{\v,\delta_{\varepsilon}\tilde{\v}-\w}^{\delta_{\varepsilon}\omega_0-\ric(\omega_0)}\) is bounded from below, then \(\M_{\v,\w}\) is \(\TC\)-coercive. 
\end{lem}

\begin{proof}
Let \(\varphi\in \mathcal{K}(X,\omega_0)^\mathbb{T}\). 
By Lemma~\ref{fut:van}, the functional \(\M_{\v,\w}\) is \(\TC\)-invariant. 
As in the proof of Theorem~\ref{t:exi:sol:delta}, let \(\gamma\in \TC\) and \(\varphi_{\gamma}\in \mathcal{K}(X,\omega_0)^\mathbb{T}\) be such that \(\omega_{\varphi_{\gamma}}=\gamma^*\omega_{\varphi}\) and 
\( (\J_{\v,\tilde{\v}}^{\omega_0})_{\TC}(\varphi) = \J_{\v,\tilde{\v}}^{\omega_0}(\varphi_{\gamma}) \). 
Then we have, by applying the Chen-Tian formula~\eqref{Chen-Tian}, the definition of \(\delta_{\varepsilon}=\delta_{\v}^{\TC}(\ka)-\varepsilon\) and the linearity property~\eqref{J-linearity} of the weighted energy functionals,  
\begin{align*}
\M_{\v,\w}(\varphi) 
& = \M_{\v,\w}(\varphi_{\gamma}) \\
& = \ent_{\v}(\varphi_{\gamma})-\J_{\v,\w}^{\ric(\omega_0)}(\varphi_{\gamma}) \\
& \geq \left(\frac{\varepsilon}{2}+\delta_{\varepsilon}\right)\J_{\v,\tilde{\v}}^{\omega_0}(\varphi_{\gamma}) -\J_{\v,\w}^{\ric(\omega_0)}(\varphi_{\gamma}) \\ 
& \geq \frac{\varepsilon}{2}(\J_{\v,\tilde{\v}}^{\omega_0})_{\TC}(\varphi) + \J_{\v,\delta_{\varepsilon}\tilde{\v}-\w}^{\delta_{\varepsilon}\omega_0-\ric(\omega_0)}(\varphi_{\gamma}) 
\end{align*}
The lemma follows. 
\end{proof}

We show a further reduction.

\begin{lem}
\label{lem_reduce_to_critical}
Assume that \(\theta\) is a Kähler form and let \(\hat{\w}\) be a smooth real valued function on \(\Delta\).  If there exists a solution \(\varphi\in \mathcal{K}(X,\omega_0)^\mathbb{T}\) to the weighted $\mathbf{J}$-equation \(\Lambda_{\varphi,\v}(\theta) = \frac{\hat{\w}(\mu_{\varphi})}{2}\), then \(\J_{\v,\hat{\w}}^{\theta}\) is bounded from below.   
\end{lem}

\begin{proof}
It was proven by Lahdili \cite[Lemma 7 and Corollary 4]{Lah20} that, if \(\theta\) is Kähler, then \(\J_{\v,\hat{\w}}^{\theta}\) is convex along weak geodesics. 
As a consequence, if \(\J_{\v,\hat{\w}}^{\theta}\) admits a critical point, then it is bounded from below. 
By Definition~\ref{def:j:v}, a solution to the weighted $\mathbf{J}$-equation \(\Lambda_{\varphi,\v}(\theta) = \frac{\hat{\w}(\mu_{\varphi})}{2}\) is a critical point. 
\end{proof}

In the remaining of the section, we will state then prove a sufficient condition of existence of a solution to the weighted $\mathbf{J}$-equation. 
Assuming Theorem~\ref{t:exi} from the next section for now, we prove our sufficient conditions of coercivity. 

\begin{proof}[Proof of Theorem~\ref{t:exi:delta}]
Let \(\hat{\w}=\delta_{\varepsilon}\tilde{\v}-\w\). 
Assumptions (i) and Lemma~\ref{lem_reduce_to_J} show that it is enough to prove that \(\J_{\v,\hat{\w}}^{\delta_{\varepsilon}\omega_0-\ric(\omega_0)}\) is bounded from below. 
By Lemma~\ref{lem_J_ddc}, it is enough to show that \(\J_{\v,\hat{\w}}^{\theta}\) is bounded from below for some \(\T\)-invariant form \(\theta\in \delta_{\varepsilon}\ka -2\pi c_1(X)\). 
By Lemma~\ref{lem_reduce_to_critical}, it suffices to show that there exists a solution \(\varphi\in \mathcal{K}(X,\omega_0)^\mathbb{T}\) to the weighted $\mathbf{J}$-equation \(\Lambda_{\varphi,\v}(\theta) = \frac{\hat{\w}(\mu_{\varphi})}{2}\) for some \(\T\)-invariant Kähler form \(\theta\in \delta_{\varepsilon}\ka -2\pi c_1(X)\).
To do so, we will apply Corollary~\ref{cor_J_eqn}, with the weights \(\v\) and \(\hat{\w}:=\delta_{\varepsilon}\tilde{\v}-\w\), and some \(\T\)-invariant Kähler form \(\theta\in \delta_{\varepsilon}\ka -2\pi c_1(X)\).

Note that the fourth assumption in Theorem~\ref{t:exi:delta} is exactly condition~\eqref{hypo:bound:2}, and the fifth assumption in Theorem~\ref{t:exi:delta} is exactly condition~\eqref{barw_convex} for \(\bar{\hat{\w}}=\check{\w}_{\varepsilon}\). 
Note that, here, \(\Delta_{\theta}\) is naturally determined, since \([\theta]=\delta_{\varepsilon}\ka-2\pi c_1(X)\) is a linear combination of \(\ka\) and \(2\pi c_1(X)\). 
To check that condition~\eqref{norma:weigh:j} holds, we do not need to fix \(\theta\in \delta_{\varepsilon}\ka-2\pi c_1(X)\) yet since it depends only on the cohomology class of \(\theta\). 
We have 
\begin{align*}
 & \int_X (2\Lambda_{0,\v}(\delta_{\varepsilon}\omega_0-\ric(\omega_0))-\hat{\w}(\mu_0))\v(\mu_0)\omega_0^{n}   \\
 & = \int_X (2\Lambda_{0,\v}(\delta_{\varepsilon}\omega_0-\ric(\omega_0))-(\delta_{\varepsilon}\tilde{\v}-\w)(\mu_0))\v(\mu_0)\omega_0^{n}   \\
& =\delta_{\varepsilon}\int_X (2\Lambda_{0,\v}(\omega_0)-\tilde{\v}(\mu_0))\v(\mu_0)\omega_0^{n} 
- 
\int_X (2\Lambda_{0,\v}(\ric(\omega_0))-\w(\mu_0))\v(\mu_0)\omega_0^{n} 
\end{align*}
The two summands being equal to zero is equivalent to the functionals \(\J_{\v,\tilde{\v}}^{\omega_0}\) and \(\J_{\v,\w}^{\ric(\omega_0)}\) being invariant under addition of a constant. 
For \(\J_{\v,\tilde{\v}}^{\omega_0}\) this follows from \cite[Lemma~3.1]{DJLb} as already recalled in Section~\ref{sec_weighted_energy}. 
For \(\J_{\v,\w}^{\ric(\omega_0)}\), this follows from condition~\eqref{norma:weight} (which is included in the vanishing weighted Futaki invariant assumption (i)) and the Chen-Tian formula. 

Now to apply Corollary~\ref{cor_J_eqn}, it remains to find \(\T\)-invariant Kähler forms \(\theta\in \delta_{\varepsilon}\ka-2\pi c_1(X)\) and \(\chi\in \ka\) such that 
\((\inf \check{\w}_{\varepsilon}) \chi - (n-1)\theta\) is Kähler.
At the level of classes, we have by using assumptions (ii) then (iii), 
\begin{align*}
(\inf \check{\w}_{\varepsilon})\ka  
& = (\inf \check{\w}_{\varepsilon}-(n-1)\delta_{\varepsilon})\ka +(n-1)\delta_{\varepsilon}\ka \\
& > (\inf \check{\w}_{\varepsilon}-(n-1)\delta_{\varepsilon})\ka +(n-1)2 \pi c_1(X) \\
& > (\inf \check{\w}_{\varepsilon}-(n-1)\delta_{\varepsilon})\ka +(n-1)\frac{-(\inf \check{\w}_{\varepsilon}-(n-1)\delta_{\varepsilon})}{n-1}\ka = 0
\end{align*} 

Since \(\ka\) is Kähler, this implies \(\inf \check{\w}_{\varepsilon}>0\). 

Let \(\theta\) be an arbitrary \(\T\)-invariant Kähler form in \(\delta_{\varepsilon}\ka-2\pi c_1(X)\). 
By assumption (iii), there exists a \(\T\)-invariant Kähler form \(\Omega\in 2\pi c_1(X)+\frac{\inf \check{\w}_{\varepsilon}-(n-1)\delta_{\varepsilon}}{n-1}\ka\). 
Consider the \(\T\)-invariant Kähler form 
\[ \chi := \frac{n-1}{\inf \check{\w}_{\varepsilon}}(\Omega+\theta) \]
We then have 
\[ (\inf \check{\w}_{\varepsilon})\chi - (n-1)\theta = (n-1)\Omega, \]
this is Kähler, and finishes the proof.  
\end{proof}

\begin{proof}[Proof of Theorem~\ref{thm_cscK_P1}]
Let \(\hat{\w}=\delta_{\varepsilon}\tilde{\v}-\w\). 
By the same first steps as in the previous proof (here, \(\delta_{\varepsilon}>1\) is equivalent to \(\theta_{\varepsilon}\) Kähler), it suffices to prove that there exists a solution to the weighted \(\J\)-equation  \(\Lambda_{\varphi,\v}(\theta) = \frac{\hat{\w}(\mu_{\varphi})}{2}\) for some \(\T\)-invariant Kähler form \(\theta\in \delta_{\varepsilon}\ka -2\pi c_1(X) = (\delta_{\varepsilon}-1)2\pi c_1(\bbP^1)\). 
Since we are in dimension one, the weighted \(\J\)-equation is actually a weighted Monge-Ampère equation  
\[ \check{\w}_{\varepsilon}(\mu_{\varphi},\mu_{\theta}) \omega_{\varphi} = \theta \]
and since we assumed \(\inf \check{\w}_{\varepsilon} > 0\) this equation admits a solution (see Section~\ref{sec_weighted_CY}). 
\end{proof}

\subsection{A sufficient condition for the weighted \(\J\)-equation}

Let $\theta$ be a $\T$-invariant K\"ahler form on $X$ with moment map \(\mu_{\theta}\) and moment image \(\Delta_{\theta}\). 
We will state here a sufficient condition of existence of a solution  \(\varphi\in \mathcal{K}(X,\omega_0)^\mathbb{T}\) to the weighted $\mathbf{J}$-equation 
\begin{equation}{\label{wt:eq}}
    \Lambda_{\varphi,\v}(\theta) = \frac{\w(\mu_\varphi)}{2}.
\end{equation}

Let $\bar{\w} : \Delta \times \Delta_\theta \rightarrow \mathbb{R}$ be the function defined by 
\begin{equation}{\label{barw}}
    \bar{\w}(x,y):= \frac{1}{2}\w(x)-\langle d\log(\v)(x),y\rangle
\end{equation} 
and recall~\eqref{hatv}.

\begin{thm}{\label{t:exi}}
There exists a solution $\varphi \in \mathcal{K}(X,\omega_0)^\T$ to  equation~\eqref{wt:eq} provided the following four conditions are satisfied:
\begin{align}
& \forall y\in \Delta_{\theta}, x\mapsto \bar{\w}(x,y) \text{ is convex} {\label{barw_convex}} \\
& \quad 1+ \inf(\hat{\v})>0 {\label{hypo:bound:2}}, \\
& \int_X (2\Lambda_{\varphi,\v}(\theta)-\w(\mu_\varphi))\v(\mu_\varphi)\omega_\varphi^{n}=0, {\label{norma:weigh:j}}
\end{align}
and  there exists a K\"ahler metric $\chi \in \ka$ such that  
\begin{align}
&\quad  \left(\inf(\bar{\w})\chi  - (n-1) \theta\right) \wedge \chi^{n-2} >0.{\label{hypo:bound}} 
\end{align}
\end{thm}

\begin{rem}
Since $\int_X \Lambda_{\varphi,\v}(\theta)\v(\mu_\varphi)\omega_\varphi^{n}$ and $\int_X (\w(\mu_\varphi))\v(\mu_\varphi)\omega_\varphi^{n}$ do not depend on $\varphi \in \mathcal{K}(X,\omega_0)^\T$ (see \cite[Lemma 2]{Lah19}), condition~\eqref{norma:weigh:j} is necessary. 
\end{rem}

\begin{rem}
In the case of cscK metrics, extremal Kähler metrics or of K\"ahler--Ricci solitons, $x \rightarrow \bar{\w}$ is convex.
\end{rem}

\begin{cor}
\label{cor_J_eqn}
Assume that the weights satisfy \eqref{barw_convex}, \eqref{hypo:bound:2} and \eqref{norma:weigh:j}. 
Assume that there exists a K\"ahler metric $\chi \in \ka$ such that 
\(\inf(\bar{\w})\chi  - (n-1) \theta\) is Kähler as well. 
Then there exists a solution $\varphi \in \mathcal{K}(X,\omega_0)^\T$ of equation~\eqref{wt:eq}.
\end{cor}

\subsection{A priori Estimates}

The estimates originate in the unweighted case from \cite{SW08} (see \cite{Wei03} for the case of surfaces) and were originally established using the $\mathbf{J}$-flow, a parabolic equation. Here, we extend an elliptic method initially proposed by \cite{LSY}.

In this section, \( C \) will denote various positive constants independent of $\varphi$, which may change from line to line.

\begin{prop}{\label{p:c2:stim}}
We suppose the same hypothesis than in Theorem \ref{t:exi}. Let $\varphi \in \mathcal{K}(X,\omega_0)^\T$ be a solution of \eqref{wt:eq}. Then there exists a uniform positive constant $C=C(\theta,\omega_0,\v,\w)$, independent of $\varphi$, such that 
\begin{equation}{\label{eq:esti}}
    \Lambda_{\theta}(\omega_\varphi) \leq C \quad \text{ and } \quad \sup\varphi-\inf \varphi \leq C.
\end{equation}
\end{prop}

By \eqref{wt:eq}, there exists $C$, $C'$ positives in $\varphi$ such that 
\begin{equation}{\label{obv:trace}}
   C \leq  \Lambda_{\varphi}(\theta) \leq \frac{1}{C} \quad \text{ and } \quad \omega_\varphi \geq C'\theta.
\end{equation}
We consider the hermitian metric $g_{\bar{\omega}_\varphi}$, which is defined for any $1$-forms  $\alpha_1$, $\alpha_2$
\begin{equation}{\label{omega:bar}}
  \quad g_{\bar{\omega}_\varphi}(\alpha_1,\alpha_2)= g_{\omega_\varphi}(\alpha_1 \wedge J\alpha_2, \theta),
\end{equation}
and we let $\bar{\omega}_\varphi$ its associate $2$-form. In particular, for any $(1,1)$-form $\beta$

\begin{equation*}
    \Lambda_{\bar{\omega}_\varphi} \beta=  g_{\omega_\varphi}(\beta, \theta). 
\end{equation*}
From \eqref{obv:trace} we also have 
\begin{equation}{\label{obv:trace2}}
    \bar{\omega}_\varphi \geq C_{\theta}\theta.
\end{equation}

We begin by the following key Lemma

\begin{lemma}{\label{l:esti:tr}}
The following inequality holds true
\begin{equation*}
\begin{split}
\Delta_{\bar{\omega}_\varphi}(\Lambda_{\theta}(\omega_\varphi)) \geq &    \frac{1}{2}\sum_{a=1}^r \w_{,a}(\mu_\varphi)\Delta_{\theta}(\mu^a_\varphi) -  \sum_{a}^r\mu_\theta^a \langle d\log(\v)_{,a}(\mu_\varphi), \Delta_\theta \mu_\varphi \rangle  \\
&+ 2|\nabla^{\theta}\omega_\varphi|^2_{\theta \otimes \bar{\omega}_\varphi } - C \Lambda_\theta(\omega_\varphi) - C.
\end{split}
\end{equation*}
\end{lemma}

\noindent In the above statement,
\begin{equation*}
|\nabla^{\theta}\omega_\varphi|^2_{\theta \otimes \bar{\omega}_\varphi }:= \sum_{j=1}^{n}|\nabla^{\theta}_{e_j}\omega_\varphi|^2_{ \bar{\omega}_\varphi } + |\nabla_{Je_j}^{\theta}\omega_\varphi|^2_{\bar{\omega}_\varphi },
\end{equation*}
where $(e_i,Je_i)_{i=1}^n$ is a normal frame for $\theta$.

\begin{proof}
From \cite[(6.3)]{LSY} (which adapt the proof of \cite[Lemma 3.1]{SW08}), we find that
\begin{equation}{\label{del:tra:theta}}
\begin{split}
   \Delta_{\bar{\omega}_\varphi}&( \Lambda_{\theta}(\omega_\varphi)) +  \Delta_{\theta}( \Lambda_{\omega_\varphi}(\theta)) \\
  =&\sum_{i,j=1}^n\frac{(|e_j|_\varphi^2-|e_i|_\varphi^2)}{|e_i|_{\bar{\omega}_\varphi}^2} \Big(R^{\theta}(e_i,e_j,e_j,e_i)+R^{\theta}(Je_i,e_j,e_j,Je_i)\Big)+2|\nabla^{\theta}\omega_\varphi|^2_{\theta \otimes \bar{\omega}_\varphi},
\end{split}  
\end{equation}
where $R^{\theta}$ denotes the curvature tensor of $\theta$ and $(e_j, Je_j)_{j=1}^n$ is a normal frame for $\theta$. Also by \cite[(3.3)]{SW08}
\begin{equation*}
 \left(\sum_{i,j=1}^n\frac{(|e_j|_\varphi^2-|e_i|_\varphi^2)}{|e_i|_{\bar{\omega}_\varphi}^2} \Big(R^{\theta}(e_i,e_j,e_j,e_i)+R^{\theta}(Je_i,e_j,e_j,Je_i)\Big)\right) \geq -C \Lambda_{\theta}(\omega_\varphi).
\end{equation*}
Hence 
\begin{equation}{\label{maj:1}}
\begin{split}
   \Delta_{\bar{\omega}_\varphi}&( \Lambda_{\theta}(\omega_\varphi)) +  \Delta_{\theta}( \Lambda_{\omega_\varphi}(\theta))
  \geq-C \Lambda_{\theta}(\omega_\varphi)+2|\nabla^{\theta}\omega_\varphi|^2_{\theta \otimes \bar{\omega}_\varphi}.
\end{split}  
\end{equation}

On the other hands, a standard computation (see e.g. \cite[(10)]{DJL}) shows that
\begin{equation*}
    \begin{split}
        \frac{1}{2}\Delta_{\theta}\left(\w(\mu_\varphi)\right) =& \frac{1}{2}\sum_{a=1}^r \w _{,a}(\mu_\varphi)\Delta_{\theta}(\mu^a_\varphi) + \frac{1}{2}\left\langle \mathrm{Hess}(\w)(\mu_\varphi), g_{\theta}(d\mu_\varphi, d\mu_\varphi) \right\rangle \\ 
    \end{split}
\end{equation*}

We compute the term$\langle d\log(\v)(\mu_\varphi), \mu_\theta \rangle= \sum_{a=1}^r \mu_\theta^a \log(\v)_{,a}(\mu_\varphi)$:

\begin{equation*}
\begin{split}
\Delta_{\theta}\big( \langle d\log(\v)(\mu_\varphi), \mu_\theta \rangle\big)= & \sum_{a=1}^r \Delta_{\theta}(\mu_\theta^a)\log(\v)_{,a}(\mu_\varphi) \\
&+ g_{\theta}( d\mu_\theta^a, d\log(\v)(\mu_\varphi))+ \mu_\theta^a \Delta_{\theta}(\log(\v)_{,a}(\mu_\varphi)).
\end{split}
\end{equation*}
 A direct computation shows that the second term of the RHS is bounded below by $-C\Lambda_{\theta}(\omega_\varphi)$. For the third term we use \cite[Eq. (10)]{DJL}
to deduce

\begin{equation*}
   \begin{split}
       -\Delta_{\theta}(\langle d\log(\v)(\mu_\varphi), \mu_\theta \rangle)\geq & -C \Lambda_{\theta}(\omega_\varphi) - \sum_{a=1}^r \mu^a_\theta \left\langle \mathrm{Hess}(\log(\v)_{,a})(\mu_\varphi), g_{\theta}(d\mu_\varphi, d \mu_\varphi)\right\rangle \\
       &-  \sum_{a}^r\mu_\theta^a \langle d\log(\v)_{,a}(\mu_\varphi), \Delta_\theta \mu_\varphi \rangle-C.
   \end{split} 
\end{equation*}
Hence since $x \rightarrow \bar{\w}(x,\cdot)$ is convex, we get that  

\begin{equation}{\label{jam}}
\begin{split}
  \frac{1}{2}\Delta_{\theta}&\left(\w(\mu_\varphi)\right) -\Delta_{\theta}(\langle d\log(\v)(\mu_\varphi), \mu_\theta \rangle)\\
  \geq&  -  \sum_{a}^r\mu_\theta^a \langle d\log(\v)_{,a}(\mu_\varphi), \Delta_\theta \mu_\varphi \rangle - C\Lambda_{\theta}(\omega_\varphi) \\
 & +\frac{1}{2}\sum_{a=1}^r \w _{,a}(\mu_\varphi)\Delta_{\theta}(\mu^a_\varphi)-C
\end{split}  
\end{equation}

We conclude using \eqref{maj:1}, \eqref{jam}:
\begin{equation*}
\begin{split}
    \Delta_{\bar{\omega}_\varphi}(\Lambda_{\theta}(\omega_\varphi)) 
    =&\Delta_{\bar{\omega}_\varphi}(\Lambda_{\theta}(\theta) + \Delta_{\theta}\left(\Lambda_{\omega_\varphi,\v}(\omega_\varphi)- \frac{1}{2}\w(\mu_\varphi) \right)\\
    \geq&  \frac{1}{2}\sum_{a=1}^r \w_{,a}(\mu_\varphi)\Delta_{\theta}(\mu^a_\varphi) -  \sum_{a}^r\mu_\theta^a \langle d\log(\v)_{,a}(\mu_\varphi), \Delta_\theta \mu_\varphi \rangle  \\
    &+ 2|\nabla^{\theta}\omega_\varphi|^2_{\theta \otimes \bar{\omega}_\varphi } - C \Lambda_\theta(\omega_\varphi) - C.
    \end{split}
\end{equation*}

\end{proof}

\begin{lemma}{\label{l:estim:log}}
There exists a positive constant $C=C(\theta,\omega_0,\v,\w)$ such that 
\begin{equation*}
\begin{split}
    \Delta_{\bar{\omega}_\varphi} \log\left( \Lambda_{\theta}(\omega_\varphi) \right) \geq &  \frac{ 1 }{\Lambda_{\theta}(\omega_\varphi)} \left( \frac{1}{2}\sum_{a=1}^r \w_{,a}(\mu_\varphi)\Delta_{\theta}(\mu^a_\varphi) -  \sum_{a}^r\mu_\theta^a \langle d\log(\v)_{,a}(\mu_\varphi), \Delta_\theta \mu_\varphi \rangle \right)-C   
\end{split}
\end{equation*}

\end{lemma}

\begin{proof}
A standard computation shows that
\begin{equation*}
    \begin{split}
\Delta_{\bar{\omega}_\varphi} \log\left( \Lambda_{\theta}(\omega_\varphi) \right)         = & \frac{ 1 }{\Lambda_{\theta}(\omega_\varphi)} \Delta_{\bar{\omega}_\varphi}  \Lambda_{\theta}(\omega_\varphi)  - \frac{|d\Lambda_{\theta}(\omega_\varphi)|^2_{\bar{\omega}_\varphi}}{\Lambda_{\theta}(\omega_\varphi)^2}. 
\end{split}
\end{equation*}

We bound the first term using Lemma \ref{l:esti:tr}:
\begin{equation*}
\begin{split}
    \frac{ \Delta_{\bar{\omega}_\varphi}  \Lambda_{\theta}(\omega_\varphi) }{\Lambda_{\theta}(\omega_\varphi)} \geq&  - \frac{C}{\Lambda_{\theta}(\omega_\varphi)} -C+\frac{ 1 }{\Lambda_{\theta}(\omega_\varphi)} 2|\nabla^{\theta}\omega_\varphi|^2_{\theta \otimes \bar{\omega}_\varphi } \\
    &+  \frac{ 1 }{\Lambda_{\theta}(\omega_\varphi)} \left( \frac{1}{2}\sum_{a=1}^r \w_{,a}(\mu_\varphi)\Delta_{\theta}(\mu^a_\varphi) -  \sum_{a}^r\mu_\theta^a \langle d\log(\v)_{,a}(\mu_\varphi), \Delta_\theta \mu_\varphi \rangle  \right) \\
   \geq  &  -C+ \frac{ 1 }{\Lambda_{\theta}(\omega_\varphi)}2|\nabla^{\theta}\omega_\varphi|^2_{\theta \otimes \bar{\omega}_\varphi } \\ 
   &+  \frac{ 1 }{\Lambda_{\theta}(\omega_\varphi)}  \left( \frac{1}{2}\sum_{a=1}^r \w_{,a}(\mu_\varphi)\Delta_{\theta}(\mu^a_\varphi) -  \sum_{a}^r\mu_\theta^a \langle d\log(\v)_{,a}(\mu_\varphi), \Delta_\theta \mu_\varphi \rangle  \right) \\
\end{split}    
\end{equation*}
where we use \eqref{obv:trace} to pass to the last line.

For the second term we use \cite[Lemma 2]{Wei03} stating that:
\begin{equation*}
   |\Lambda_{\theta}(\nabla^{\theta}\omega_\varphi)|^2_{\bar{\omega}_\varphi} \leq  2\Lambda_{\theta}(\omega_\varphi)|\nabla^{\theta}\omega_\varphi|^2_{\theta \otimes \bar{\omega}_\varphi},
\end{equation*}
which concludes the proof.

\end{proof}

\begin{proof}(of Proposition \ref{p:c2:stim}).
Let $A>0$ be a positive constant. Consider a normal frame for $\theta$. Using that $\varphi$ solves \eqref{wt:eq} we find that
\begin{equation}{\label{comp:delta}}
\begin{split}
  -A\Delta_{\bar{\omega}_\varphi}\varphi
  =&  -A\Lambda_{\bar{\omega}_\varphi}(\omega_\varphi-\omega_0) \\
      &= -A\Lambda_{\omega_\varphi}(\theta)+A\Lambda_{\bar{\omega}_\varphi}(\omega_0)  \\
     &= -A\bar{\w}(\mu_\varphi,\mu_\theta) +A\Lambda_{\bar{\omega}_\varphi}(\omega_0)  .
\end{split}      
\end{equation}

We will apply the maximum principle. Let 
   $ x_0 \in X$ be a maximal point of  $\log\left( \Lambda_{\theta}(\omega_\varphi) \right) -  A \varphi)$. At $x_0$,  for any $\xi_a$ we have that
\begin{equation*}
   \mathcal{L}_{J\xi_a} \left( \log\left( \Lambda_{\theta}(\omega_\varphi) \right) - A \varphi \right)=0.
\end{equation*}

Developing the LHS of the above equality, we get
\begin{equation*}
    \begin{split}
   \mathcal{L}_{J\xi_a} \left( \log\left( \Lambda_{\theta}(\omega_\varphi) \right) - A \varphi \right) =&\frac{\nabla^\theta_{J\xi_a}\left( \Lambda_{\theta}(\omega_\varphi) \right)}{\Lambda_{\theta}(\omega_\varphi)}  +A d^c\varphi(\xi_a)  \\
  = &  \frac{ \Lambda_{\theta}(\nabla^\theta_{J\xi_a}\omega_\varphi) }{\Lambda_{\theta}(\omega_\varphi)}   +  Ad^c\varphi(\xi_a)\\    
=& - \frac{\Delta_{\theta}\mu_\varphi^a }{\Lambda_{\theta}(\omega_\varphi)}+ \frac{1}{\Lambda_{\theta}(\omega_\varphi)}\left(  \Lambda_\theta\left( g_{\omega_\varphi} ( \nabla^\theta\xi_a, \cdot) \right)^{\rm skw}\right) \\
    &+ A (\mu_\varphi^a-\mu_0) \\
\end{split}
\end{equation*}
The first equality is due to the fact that the Lie derivative and the covariant derivative coincide on function. For the second equality, we use \cite[Lemma 5.4]{DJL}. For the third one, we use that $ \nabla^{\theta}_{J\xi_a}(\omega_\varphi)= -dd^c\mu^a_\varphi + 2 \left( g_{\omega_\varphi} ( \nabla^\theta\xi_a, \cdot) \right)^{\rm skw}$ (see \cite[p.17]{DJL}). Using  that $C \Lambda_\theta(\omega_\varphi) \geq \Lambda_{\theta}\left(  g_{\omega_\varphi} ( \nabla^\theta\xi_a, \cdot) \right)^{\rm skw} \geq -C \Lambda_\theta(\omega_\varphi)$ (see \cite[p.18]{DJL}) and \eqref{obv:trace}, we deduce that at $x_0$
\begin{equation}{\label{lap:moment}}
\begin{split}
-  \frac{|\Delta_{\theta}\mu^a_\varphi |}{\Lambda_{\theta}(\omega_\varphi)} 
     \geq - C.
\end{split}     
\end{equation}

From Lemma \ref{l:estim:log}, \eqref{comp:delta} and applying the maximum principle  we deduce that at $x_0$ 
\begin{equation*}
\begin{split}
    0 \geq   \Delta_{\bar{\omega}_\varphi}(\log(\Lambda_{\theta}(\omega_\varphi)-A\varphi+ \bar{\v}(\mu_\varphi))
    \geq   -C -A \bar{\w}(\mu_\varphi,\mu_\theta) +A \Lambda_{\bar{\omega}_\varphi}(\omega_0).
\end{split}    
\end{equation*}
Hence writing $\epsilon:= \frac{C}{A}$, we obtain
\begin{equation*}
\epsilon \geq  - \bar{\w}(\mu_\varphi,\mu_\theta) +  \Lambda_{\bar{\omega}_\varphi}(\omega_0).
\end{equation*}

Now the proof is similar than in the cscK case (see \cite[Lemma 3.1]{SW08}), but we give some explanations for reader's convenience. Let $(e_i,Je_i)_{i=1}^n$ be a normal frame for $\omega_0$.  By \cite[(6.6))]{LSY} for any index $k=1,\dots, n$,
\begin{equation}{\label{ineq:eps}}
    \epsilon  \geq \inf_{\Delta \times \Delta_{\mu_\theta}}   (\bar{\w})  - \sum_{i\neq k} |e_i|^2_{\theta} - 2\frac{|e_k|^2_{\theta}}{|e_k|^2_{\varphi}}.
\end{equation}

Applying our hypothesis \eqref{hypo:bound} and choosing $\omega_0=\chi$, we find that 
\begin{equation*}
\left(\inf_{\Delta \times \Delta_{\mu_\theta}}(\bar{\w})\omega_0   - (n-1) \theta  \right)\wedge\omega_0^{n-2}\wedge( e_k \wedge J e_k) > 3\epsilon \omega_0^{n-1} \wedge ( e_k \wedge J e_k),
\end{equation*}
for $\epsilon$ small enough. Hence by \cite[(6.9)]{LSY} (see also \cite{SW08})
\begin{equation*}
\inf_{\Delta \times \Delta_{\mu_\theta}}   (\bar{\w})  - \sum_{i\neq k} |e_i|^2_{\theta} >  3\epsilon
\end{equation*}
Substituting back in \eqref{ineq:eps} we get that at the maximal point $x_0$
\begin{equation*}
  \frac{|e_k|^2_{\varphi}}{|e_k|^2_{\theta}} \leq \frac{1}{\epsilon},
\end{equation*}
for every $k=1,\dots, n$. Then at $x_0$
\begin{equation*}
  \Lambda_\theta(\omega_\varphi) \leq \frac{n}{\epsilon}.
\end{equation*}
Hence, we deduce the $C^2$-estimates on  $X$ 
\begin{equation}{\label{c2:est}}
    \Lambda_\theta(\omega_\varphi) \leq \frac{n}{\epsilon}e^{A(\varphi-\inf \varphi)}.
\end{equation}

The demonstration of the $C^0$ from the $C^2$-estimates do not use the equation \eqref{wt:eq} but only \eqref{c2:est}. Hence we can apply the same arguments than \cite{Wei03} and \cite{SW08}.
\end{proof}

\subsection{Continuity method}
The idea of employing a continuity path combined with an elliptic operator is introduced in \cite[Theorem 1.6]{LSY} within the context of the cscK case. This approach forms the basis of the generalization presented in this section.

We consider the weighted Li--Shi--Yao continuity path 
\begin{equation}{\label{cont:path2}}
    \Lambda_{\varphi,\v}(\theta_t)=\frac{1}{2}\w_t(\mu_\varphi), \quad \varphi \in \mathcal{K}(X,\omega_0)^\T, \text{ } t \in [0,1],
\end{equation}
where $\theta_t := t \theta + (1-t) \omega_0$ and $\w_t(x)=t \w(x) + (1-t)\tilde{\v}(x) $, where $\tilde{\v}$ is introduced in \eqref{defn_tildev}. Clearly a solution at $t=1$ satisfies \eqref{wt:eq}. We consider
\begin{equation*}
    S:= \{ t \in [0,1] \text{ } | \text{ } \exists    \varphi_t \in \mathcal{K}(X,\omega_0)^\T \text{ solution of } \eqref{cont:path2} \}
\end{equation*}

The set $S$ is non-empty since $\varphi=0$ is a solution at $t=0$ of \eqref{cont:path2}. We prove that $S$ is open and closed.

First we observe that for any $t\in[0,1]$,
\begin{equation}{\label{norma:cont:path}}
\begin{split}
        \int_X&\left(2\Lambda_{\varphi,\v}(\theta_t)-\w_t(\mu_{\varphi})\right)\v(\mu_{\varphi})\omega_{\varphi}^{[n]} \\
        =& t  \int_X\left(2\Lambda_{\varphi,\v}(\theta)-\w(\mu_{\varphi})\right)\v(\mu_{\varphi})\omega_{\varphi}^{[n]}   \\
        &+ (1-t)  \int_X\left(2\Lambda_{\varphi,\v}(\omega_0)-\tilde{\v}(\mu_{\varphi})\right)\v(\mu_{\varphi})\omega_{\varphi}^{[n]} \\
        =&0.
\end{split}        
\end{equation}
The first term of the right-hand side is zero by hypothesis \eqref{norma:weigh:j}, and second term of the right-hand side is zero since the integral does not depend on $\varphi_t \in \mathcal{K}(X,\omega_0)^\T$, see \cite[Lemma 2]{Lah19}.

We consider the operator $  R : \mathcal{K}(X,\omega_0)^\T \times [0,1] \longrightarrow \mathcal{C}^\infty(X)^\T$ defined by
\begin{equation*}
    R(\varphi,t):= \Lambda_{\varphi,\v}(\theta_t)-\frac{1}{2}\w_t(\mu_\varphi),
\end{equation*}
Observe that by \eqref{norma:cont:path}, $R(\varphi,t)$ is $L^2(\v(\mu_{\varphi})\omega_{\varphi}^{[n]})$-orthogonal to constants in $\mathcal{C}^\infty(X)^\T$.

We fix $(\varphi_0,t_0)$ a solution of (\ref{cont:path2}) and we let $\omega_{t_0}:=\omega_0+dd^c\varphi_{0}$. 

\begin{lemma}{\label{l:line}}
 The linearization of $R$ in the direction $(\dot{\varphi},s)$ at $(\varphi_{t_0},t_0)$ is given by 
\begin{equation}{\label{hash}}
D_{(\varphi_{0},t_0)}(R)[\dot{\varphi},s]=\mathbb{H}^{\theta_{t_0}}_{\varphi_0,\v}(\dot{\varphi}) + s\big(\Lambda_{\varphi_0,\v}(\theta)- \w(\mu_{\varphi_0}) - (\Lambda_{\varphi_0,\v}(\omega_0)-\tilde{\v}(\mu_{\varphi_0}))\big),
\end{equation}
where
\begin{equation*}
  \mathbb{H}^{\theta}_{\varphi,\v}f:= g_{\omega_\varphi}\big(\theta, dd^c f \big) + g_{\omega_\varphi}\big(d \Lambda_{\varphi} \theta, df\big)  + g_{\omega_\varphi}\big(\theta,d\log(\v(\mu_{\varphi})) \wedge d^cf\big).
\end{equation*}
In particular, if $\theta$, $\v$, $\w$ satisfy \eqref{norma:weigh:j}, $S$ is open.
 \end{lemma}

\begin{proof}
Equation \eqref{hash} is proved along the lines of the proof of \cite[Lemma 3.5]{DJLb}. By \cite[Lemma 3.4]{DJLb}, $\mathbb{H}^{\theta}_{\varphi,\v}$ is elliptic and self-adjoint with respect to $\v(\mu_\varphi)\omega_\varphi^{[n]}$,  hence the conclusion follows from standard arguments of elliptic theory. 
\end{proof}

\begin{lemma}
Assume that $\theta$, $\v$, $\w$ satisfy the hypothesis of Theorem \ref{t:exi}, and set $\omega_0=\chi$, where \(\chi\) is given by assumption~\eqref{hypo:bound}. Then $S$ is closed.
\end{lemma}

\begin{proof}
We need to show that the estimates~\eqref{eq:esti} hold along the continuity path. The proof is essentially the same than Proposition~\ref{p:c2:stim} (see \cite[Proof of Lemma~6.1]{LSY} for a detailed proof in the cscK case) so we do not reproduce it here. We only check that $\theta_t$, $\v$, $\w_t$ satisfy condition \eqref{hypo:bound} for all \(t\). 

 Let $\mu_t:=\mu_{\varphi_t}$. We compute
\begin{equation*}
    \begin{split}
        \frac{1}{2}\w_t(\mu_t)- &\langle d \log(\v)(\mu_t), \mu_{\theta_t} \rangle \\
        &= t\left( \frac{1}{2} \w(\mu_t) -  \langle d \log(\v)(\mu_t), \mu_{\theta} \rangle \right) + (1-t)\left( \tilde{\v}(\mu_t) -  \langle d \log(\v)(\mu_t), \mu_{0} \rangle \right)
    \end{split}
\end{equation*}
Then 
\begin{equation}{}\label{ineq:closed}
\begin{split}
    & \left(\left(\frac{1}{2}\w_t(\mu_t)- \langle d \log(\v)(\mu_t), \mu_{\theta_t} \rangle\right) \wedge \omega_0 - (n-1)\theta_t)\right)\wedge \omega_0^{n-2} \\
    &= t\left(\left(\frac{1}{2}\w(\mu_t)- \langle d \log(\v)(\mu_t), \mu_{\theta} \rangle\right) \wedge \omega_0 - (n-1)\theta)\right)\wedge \omega_0^{n-2}\\
    &+ (1-t)\left(\left(\frac{1}{2}\tilde{\v}(\mu_t)- \langle d \log(\v)(\mu_t), \mu_{0} \rangle\right) \wedge \omega_0 - (n-1)\omega_0)\right)\wedge \omega_0^{n-2}
\end{split}
\end{equation}

For any $t>0$, the first term of the RHS is positive by hypothesis. For the second term, we have that
\begin{equation*}
    \begin{split}
   (1-t)&\left(\left(\frac{1}{2}\tilde{\v}(\mu_t)- \langle d \log(\v)(\mu_t), \mu_{0} \rangle\right) \wedge \omega_0 - (n-1)\omega_0)\right)\wedge \omega_0^{n-2} \\
    &=(1-t)\left(1 + \langle d \log(\v)(\mu_t), d^c\varphi_t \rangle\right) \omega_0^{n-1},
    \end{split}
\end{equation*}
which is positive by \eqref{hypo:bound:2}. Now if $t=0$, the above equality, together with \eqref{hypo:bound:2} show that the RHS of \eqref{ineq:closed} is positive.

\end{proof}

\section{weighted beta invariant for a semisimple principal fibration}
\label{sec_fibrations}

\subsection{Setting}
\label{sec_sspf}

Let \(B_a\) be a finite collection of projective smooth manifolds, with \(a\in \{1,\ldots,k\}\). 
For each \(a\), fix an ample line bundle \(L_a\), and denote by \(L_a^0\) its zero section. 
Then the product \(Q:= \prod_a L_a \setminus L_a^0\) defines a principal \((\bbC^*)^k\)-bundle over the base \(\prod_a B_a\). 
Fix an integer \(r\in \bbZ_{\geq 1}\) and for each \(a\), a one-parameter subgroup \(p_a\) of \((\bbC^*)^r\). 
Finally, let \(X\) be a projective smooth manifold, equipped with a holomorphic action of \(\T=(\bbC^*)^r\). 

We call \emph{semisimple principal bundle} associated to this data the fiber bundle 
\[ Y=(X\times Q)/(\bbC^*)^k \] 
where \((\bbC^*)^k\) acts diagonally, with its natural action on \(Q\), and through the given \((\bbC^*)^r\)-action on \(X\) via the morphim \((p_a)_{a=1}^k:(\bbC^*)^k\to (\bbC^*)^r\). 
It is a complex manifold equipped with a holomorphic projection map \(\pi=(\pi_a):Y\to B=\prod_a B_a\) and a holomorphic fiberwise action of \(\T\). 

Choose positively curved Hermitian metrics \(h_a\) on each \(L_a\), and corresponding Kähler curvature forms \(\omega_a\in 2\pi c_1(L_a)\). 
Note that by considering the product of \(\bbS^1\)-bundles defined by \(h_a=1\) in \(L_a\), we obtain an underlying \((\bbS^1)^k\)-principal bundle \(P\) in \(Q\) which defines a restriction of the structure group. 
The manifold \(Y\) may as well be written as the quotient \(X\times P /(\bbS^1)^k\) under the action defined as above by the inclusion \((\bbS^1)^k\subset (\bbC^*)^k\). 

The metrics also naturally define a connection \(\theta\) on \(Q\), which descends to \(Y\). 
Given a \(\T=(\bbS^1)^r\)-invariant real closed \((1,1)\)-form \(\omega\) on \(X\), and a fixed (\(\T\)-invariant) moment map \(\mu_{\omega}\), 
 we define a global closed real \((1,1)\)-form on \(Y\) by 
\begin{equation}
\label{eqn_induced_form}
\omega + \sum_{a=1}^k\langle p_a, \mu_\omega \rangle \pi_{B_a}^*(\omega_a) 
\end{equation}
where by abuse of notations we mean that if \(\xi^{(i)} = \xi^{(i)}_h + \xi^{(i)}_v\) are tangent vector fields at the same point \([x,p]\), decomposed into horizontal and vertical part with respect to the connection, the term \(\omega\) in the above \((1,1)\)-form applied to \((\xi^{(1)},\xi^{(2)})\) reads
\( \omega_x(\xi_v^{(1)},\xi_v^{(2)}) \) and the moment map \(\mu_{\omega}\) is evaluated at \(x\). 
This is well-defined by invariance with respect to \(\T\), we refer to \cite{ACGT_2004,ACGT_2011,AJL} for more details on this construction. 
By a further abuse of notation, we will omit from the notation the pullback in \(\pi_{B_a}^*\omega_a\) from now on. 

\begin{rem}
We consider in this section semisimple principal fibration in a weaker sense than in previous papers (e.g. \cite{AJL, Jubert_2023,Delcroix-Jubert_2023}): we do not assume that the basis is equipped with a cscK metric. 
Furthermore, we assume throughout that \(Y\) has a projective fiber in order to highlight easily the various possible choices of connections when a semisimple principal fibration is given a fixed complex structure. 
As a consequence, we provide a definition of compatible Kähler metrics which is somewhat different from previous definitions. 
\end{rem}

\begin{defn}
\label{defn_comp}
A \(\T\)-invariant Kähler form \(\tilde{\omega}\) on \(Y\) is called \emph{compatible} if it writes as 
\[ \tilde{\omega} = \omega + \sum_{a=1}^k(\langle p_a, \mu_\omega \rangle+c_a)\omega_a \]
for some choice of: positively curved Hermitian metrics \(h_a\) with curvature form \(\omega_a\), \(\T\)-invariant Kähler form \(\omega\) on \(X\) with moment map \(\mu_{\omega}\), and constants \(c_a\). 

We say that a Kähler class \([\tilde{\omega}_0]\) on \(Y\) is \emph{compatible} if it contains a compatible Kähler metric. 
\end{defn}

Note that, if \(\tilde{\omega}\) is Kähler, then \(\langle p_a, \mu_\omega \rangle+c_a\) must be a positive function on \(X\). Conversely, if it is, then the associated closed real \((1,1)\)-form is Kähler.  

The Ricci form of a compatible Kähler metric may be written similarly as the compatible metric itself, thanks to the expression below proved in \cite{AJL} as a step in the proof of Lemma~5.9. 
Note that, at this stage in the proof of \cite[Lemma~5.9]{AJL}, the assumption that the basis factors are cscK has not been used. 

\begin{prop}[{\cite[Equation~(40),p.3264]{AJL}}]{\label{p:ricci}}
The Ricci form of a compatible K\"ahler metric $\tilde{\omega}$ is given by
\begin{equation*}
    \ric(\tilde{\omega}) = \ric_\p(\omega) + \sum_{a=1}^k\left(\ric(\omega_a) + \langle p_a, \mu_{\ric_{\p}(\omega) }\rangle \omega_a\right),
\end{equation*}
where \(\p:\mu_{\omega}(X)\to \bbR, y\mapsto \prod_a(\langle p_a,y\rangle +c_a)^{\dim B_a}\) and \(\mu_{\ric_\p(\omega)}\) is the moment map of the closed real \((1,1)\)-form \(\ric_\p(\omega)\in 2\pi c_1(X)\). 
\end{prop}

By the definition of the weighted Ricci curvature and \eqref{eqn_induced_form}, we deduce the expression of the weighted Ricci curvature as well. 

\begin{cor}
Given a weight \(\v:\mu_{\omega}(X)\to \bbR_{>0}\), we have 
\begin{equation*}
    \ric_{\v}(\tilde{\omega}) = \ric_{\v\p}(\omega) + \sum_{a}\ric(\omega_a) + \langle p_a, \mu_{\ric_{\p\v}(\omega) }\rangle \omega_a
\end{equation*}
\end{cor}

\subsection{Greatest compatible weighted Ricci lower bounds of fibrations}

In this section, we consider greatest weighted Ricci lower bound for both the fiber \(X\) and the fiber bundle \(Y\), so to make things more explicit, we will include the manifold in the notation. 
For example, \(\beta(Y,[\tilde{\omega}_0])\) denotes the greatest Ricci lower bound for the Kähler class \([\tilde{\omega}_0]\) on \(Y\). 
Furthermore, we introduce a variant of the greatest weighted Ricci lower bound for compatible Kähler classes on a semisimple principal fiber bundle. 

\begin{defn}
If \(Y\) is a semisimple principal fibration as before, and \([\tilde{\omega}_0]\) is a compatible Kähler class, we set 
\[ \beta^{\compatible}_{\v}(Y,[\tilde{\omega}_0]) = \sup \{\beta \in \bbR\mid \exists 0 < \tilde{\omega} \in [\tilde{\omega}_0] \text{ compatible st } \ric_{\v}(\tilde{\omega})-\beta\tilde{\omega}>0 \} \]
\end{defn}

In the next statement, we compute \(\beta^{\compatible}_{\v}(Y,[\tilde{\omega}_0])\) in terms of the basis and the weighted fiber. 
The general expression of \(\beta^{\compatible}_{\v}(Y,[\tilde{\omega}_0])\) is somewhat involved, but can be drastically simplified once the fiber is fixed and say Fano, as we will see in the next section.  

\begin{thm}
\label{thm:beta_comp}
With the notations as above, we have 
\[ \beta^{\compatible}_{\v}(Y,[\tilde{\omega}_0]) = \inf_a \sup\{ t<\beta_{\v\p}(X,\ka) \mid \inf_{\Delta_t} \langle p_a ,\cdot \rangle -tc_a > - \beta(B_a,[\omega_a]) \} \] 
where \(\Delta_t=\Delta_{2\pi c_1(X)-t\ka}\) is the moment image of the Kähler class \(2\pi c_1(X)-t\ka\) 
\end{thm}

\begin{proof}
We decompose, for \(\beta\in \bbR\), the form \(\ric_{\v}(\tilde{\omega})-\beta\tilde{\omega}\) along the vertical directions and the horizontal directions, the latter being further decomposed along the product structure on the basis. 
We obtain the expression 
\[ (\ric_{\v\p}(\omega) - \beta \omega) + \sum_{a=1}^k\left(\ric(\omega_a) + (\langle p_a, \mu_{\ric_{\p\v}(\omega) }\rangle -\beta(\langle p_a,\mu_{\omega}\rangle +c_a))\omega_a\right)\]
We see readily from this expression (recall our abuse of notations), that this form is positive if and only 
\begin{enumerate}
\item \(\ric_{\v\p}(\omega) - \beta \omega\) is positive on \(X\), 
\end{enumerate}
and for each \(x\in X\) and each \(a\),
\begin{enumerate}[resume]
\item \(\ric(\omega_a) + (\langle p_a, \mu_{\ric_{\p\v}(\omega)}(x)-\beta\mu_{\omega}(x)\rangle -\beta c_a)\omega_a\) is positive on \(B_a\). 
\end{enumerate}
By the first condition, if \(\ric_{\v}(\tilde{\omega})-\beta\tilde{\omega}\) is Kähler, then \(\ric_{\v\p}(\omega) - \beta \omega\) is Kähler on \(X\) and its moment image \(\Delta_t\) is not empty. 
Now the second condition translates as the expression in the statement of the theorem since 
\[ \Delta_t = \{\mu_{\ric_{\p\v}(\omega)}(x)-\beta\mu_{\omega}(x)\mid x\in X\} \]
\end{proof}

\subsection{The anticanonical class and compatibly Fano fibrations}

Let us consider the case, of special interest, of Fano fibrations. 
More precisely: 

\begin{defn}
A semisimple principal fibration \(Y\) is called \emph{compatibly Fano} if there exists a compatible Kähler metric in \(2\pi c_1(Y)\) such that \(\ric(\tilde{\omega})\) is Kähler. 
\end{defn}

From Proposition~\ref{p:ricci}, we have the following characterization. 

\begin{prop}
The semisimple principal fibration \(Y\) is compatibly Fano if and only if \(X\) is Fano, and for all \(a\), \(c_a [\omega_a]= 2\pi c_1(B_a)\) with \(c_a>0\) (hence \(B_a\) is Fano), and \(\inf_{\Delta_{2\pi c_1(X)}}\langle p_a,\cdot\rangle > -c_a \beta(B_a,c_1(B_a))\), where \(\Delta_{2\pi c_1(X)}\) denotes the moment image of \(2\pi c_1(X)\). 
\end{prop}

\begin{proof}
Assume that there exists a compatible Kähler metric \(\tilde{\omega}\) in \(2\pi c_1(Y)\) such that \(\ric(\tilde{\omega})\) is Kähler. 
Then by Proposition~\ref{p:ricci}, the form \(\ric_{\p}(\omega)\) is Kähler on \(X\), and for each \(x\in X\), and each \(a\), the form \(\ric(\omega_a)+\langle p_a,\mu_{\ric_{\p}(\omega)}(x)\rangle \omega_a\) is Kähler on \(B_a\). 
Thus we have that \(X\) is Fano, and that \(\inf_{\Delta_{2\pi c_1(X)}}\langle p_a,\cdot,\rangle > -\beta(B_a,[\omega_a])\) for each \(a\). 
Furthermore, we have that the classes of \(\tilde{\omega}\) and of \(\ric(\tilde{\omega})\) are both equal to \(2\pi c_1(Y)\), which by the decompositions yield 
\( [\omega] = [\ric_{\v}(\omega)] =2\pi c_1(X)\) (and we can assume the moment image of \([\omega]\) is \(\Delta_{2\pi c_1(X)}\), as that of \(\ric_\v(\omega)\)) and \(c_a [\omega_a] = [\ric(\omega_a)] = 2\pi c_1(B_a)\). 
By the Kähler condition, we have that \(\langle p_a,\cdot \rangle + c_a\) is positive on \(\Delta_{2\pi c_1(X)}\). 
Since \(\Delta_{2\pi c_1(X)}\) always contains the origin, \(c_a>0\), thus \(B_a\) is Fano as well. 
Note that the condition \(\inf_{\Delta_{2\pi c_1(X)}}\langle p_a,\cdot,\rangle > -\beta(B_a,[\omega_a])\) is stronger than the condition \(\inf_{\Delta_{2\pi c_1(X)}}\langle p_a,\cdot,\rangle > -c_a\) since 
\begin{equation} 
\label{eqn:basis_Fano}
\beta(B_a,[\omega_a]) = \beta(B_a,\frac{
1}{c_a}2\pi c_1(B_a)) = c_a \beta(B_a,2\pi c_1(B_a)) \leq c_a 
\end{equation}

The converse is straightforward. 
\end{proof}

\begin{rem}
The last condition in the characterization above is a subtlety introduced by the requirement that there is a compatible Kähler form with positive Ricci form in \(2\pi c_1(Y)\). 
It is stronger than simply requiring \(Y\) to be Fano as we prove by exhibiting Example~\ref{exa_not_comp_Fano}. 
It is a striking illustration of the relevance of the greatest Ricci lower bound for dealing with fibrations. 
On the other hand, if \(B\) is Kähler-Einstein, then it follows from \cite[Lemma~5.11]{AJL} that \(Y\) is compatibly Fano if and only if it is Fano.  
\end{rem}

Our main statement yields the following simplified expression in the compatibly Fano case, and more generally, when we merely assume that the restriction to the fibers is a multiple of the anticanonical class (so that \(X\) is Fano, but \(Y\) may not be Fano at all). 

\begin{prop}
\label{prop:FanoFiber}
Assume that \(\ka = \lambda 2\pi c_1(X)\), \(\lambda>0\), with moment image \(\lambda \Delta_{2\pi c_1(X)}\). Then 
\[ \beta_{\v}^{\compatible}(Y,[\tilde{\omega}_0]) = \min \left\{ \beta_{\v\p}(X,\ka), \min_a \frac{\beta(B_a,[\omega_a])+\inf_{\Delta_{2\pi c_1(X)}}p_a}{c_a+\lambda\inf_{\Delta_{2\pi c_1(X)}}p_a} \right\} \]
In particular, if \(Y\) is compatibly Fano, then 
\[ \beta_{\v}^{\compatible}(Y,2\pi c_1(Y)) = \min \left\{ \beta_{\v\p}(X,2\pi c_1(X)), \min_a \frac{c_a \beta(B_a,2\pi c_1(B_a))+\inf_{\Delta_{2\pi c_1(X)}}p_a}{c_a+\inf_{\Delta_{2\pi c_1(X)}}p_a} \right\} \] 
\end{prop}

\begin{proof}
Observe that the moment image of \(2\pi c_1(X)-t\ka = (1-t\lambda)2\pi c_1(X)\) is \((1-t\lambda)\Delta_{2\pi c_1(X)}\). 
In the expression given by Theorem~\ref{thm:beta_comp}, we can thus replace \(\min_{\Delta_t} p_a\) with \((1-t\lambda)\min_{\Delta_{2\pi c_1(X)}} p_a\). 
This allows us to isolate the \(t\) in that expression and get: 
\[ \beta_{\v}^{\compatible}(Y,[\tilde{\omega}_0]) = \min_a \sup\left\{ t<\beta_{\v\p}(X,\ka) \mid t < \frac{\beta(B_a,[\omega_a])+\inf_{\Delta_{2\pi c_1(X)}}p_a}{c_a + \lambda \inf_{\Delta_{2\pi c_1(X)}}p_a} \right\} \]
hence the result.
If \(Y\) is furthermore compatibly Fano, then \(c_a[\omega]=2\pi c_1(B_a)\), and equation~\eqref{eqn:basis_Fano} gives the further simplification. 
\end{proof}

\begin{cor}
If \(Y\) is compatibly Fano and \(\beta_{\v}^{\compatible}(Y,2\pi c_1(Y)) = 1\), then 
\[ \beta_{\v}(Y,2\pi c_1(Y))
= \beta_{\v\p}(X,2\pi c_1(X))
= \beta(B_a,2\pi c_1(B_a))=1. \] 
\end{cor}

\begin{proof}
This follows directly from the expression above, together with the fact that \(c_a>0\) and \(\beta(B_a,2\pi c_1(B_a))\leq 1\). 
\end{proof}

\begin{cor}
\label{cor:K-ss_basis}
If \(Y\) is compatibly Fano and the \(B_a\) are K-semistable, then 
\[ \beta_{\v}^{\compatible}(Y,2\pi c_1(Y)) = \beta_{\v\p}(X,2\pi c_1(X)). \]
\end{cor}

\begin{proof}
Under the assumption, we have 
\(\beta(B_a,2\pi c_1(B_a))=1\) hence 
\[ \frac{c_a \beta(B_a,2\pi c_1(B_a))+\inf_{\Delta_{2\pi c_1(X)}}p_a}{c_a+\inf_{\Delta_{2\pi c_1(X)}}p_a}=1. \] 
Since \(\beta_{\p\v}(X,2\pi c_1(X))\leq s(X,2\pi c_1(X))=1\) we get the result. 
\end{proof}

If the fiber is toric, our general upper bound for the weighted beta invariant shows that our lower bound may easily become sharp and actually computes the weighted beta invariant rather than only the compatible weighted beta invariant. 

\begin{prop}
\label{prop_Fano_toric_fiber}
If \(Y\) is compatibly Fano, \(X\) is toric and 
\[ \beta_{\v\p}(X,2\pi c_1(X)) \leq \min_a \frac{c_a \beta(B_a,2\pi c_1(B_a))+\inf_{\Delta_{2\pi c_1(X)}}p_a}{c_a+\inf_{\Delta_{2\pi c_1(X)}}p_a}, \] 
then 
\[ \beta_{\v}^{\compatible}(Y,2\pi c_1(Y)) = \beta_{\v}(Y,2\pi c_1(Y)) = \beta_{\v\p}(X,2\pi c_1(X)). \]
\end{prop}

\begin{proof}
Under the assumptions, using Proposition~\ref{prop:FanoFiber}, we have 
\[ \beta_{\v\p}(X,2\pi c_1(X)) = \beta_{\v}^{\compatible}(Y,2\pi c_1(Y)) \]
We obviously have 
\[ \beta_{\v}^{\compatible}(Y,2\pi c_1(Y)) \leq \beta_{\v}(Y,2\pi c_1(Y)) \]
by definition. 
Finally, the upper bound from Theorem~\ref{thm_general_upper} gives, under our assumptions, 
\[ \beta_{\v}(Y,2\pi c_1(Y)) \leq \beta_{\v\p}(X,2\pi c_1(X)). \]
Putting together this sequence of inequalities, we obtain the equality in the statement. 
\end{proof}

\subsection{On simple principal \(\bbP^1\)-bundles}

\subsubsection{Setting} 

To illustrate how our main statement may easily be applied explicitly in examples, we consider the case when \(k=1\) and \(X=\bbP^1\), equipped with its toric structure under the action of \(\bbC^*\). 
We identify  one-parameter subgroups of \(\bbC^*\) with elements of \(\bbZ\), so that \(p:=p_1\) is an integer (we omit the index \(1\) for simplicity). 
Furthermore the moment images of Kähler classes on \(\bbP^1\) lie in a one-dimensional vector space identified accordingly with \(\bbR\). 
We can write \(\ka = \lambda 2\pi c_1(\bbP^1)\) for some \(\lambda >0\) and without loss of generality, we assume that the moment image of \(\ka\) is \([-\lambda,\lambda]\). 
By the compatible Kähler property, we have, for any  \(u\in [-\lambda,\lambda]\),  \(pu+c>0\), or in other words, \(c>\lambda \lvert p\rvert\).  
Furthermore, \(\p(u)=(pu+c)^{d}\) where \(d=\dim(B)\). 
By symmetry, we may as well assume that \(p\geq 0\), and since the case \(p=0\) corresponds to products, we assume \(p>0\). 

By Proposition~\ref{prop:FanoFiber}, we have

\[ \beta^{\compatible}(Y,[\tilde{\omega}_0]) = \min\left\{\beta_{\p}(\bbP^1,\ka), \frac{\beta(B_1,[\omega_1])-\lvert p\rvert}{c-\lambda\lvert p\rvert} \right\} \]

\subsubsection{The weighted \(\bbP^1\)}

We now compute the quantity \(\beta_{\p}(\bbP^1,\ka)\). 
By equation~\eqref{eqn:beta_multiple}, we have 
\[ \beta_{\p}(\bbP^1,\ka) = \beta_{\p}(\bbP^1,\lambda 2\pi c_1(\bbP^1)) = \frac{1}{\lambda}\beta_{\p(\lambda \cdot)}(\bbP^1,2\pi c_1(\bbP^1)) \]
The latter can be computed effectively 
(the formula for arbitrary log concave weights must have been first proven in \cite[Theorem~1.2]{Delcroix-Hultgren_2021}):
\begin{align*}
\beta_{\p}(\bbP^1,\ka) 
& = \frac{1}{\lambda}\sup\left\{t\in [0,1[ ; \frac{-t}{(1-t)}\bary_{\p(\lambda\cdot)}([-1,1]) \in [-1,1]\right\} \\
& = \frac{1}{\lambda}\sup\left\{t\in [0,1[ ; \left\lvert \frac{-t}{(1-t)}\bary_{\p(\lambda\cdot)}([-1,1]) \right\rvert \leq 1 \right\} \\
& = \frac{1}{\lambda}\frac{1}{1+\lvert \bary_{\p(\lambda\cdot)([-1,1])} \rvert} 
\end{align*} 
where 
\begin{align*} 
\bary_{\p(\lambda\cdot)}([-1,1]) & = \frac{\int_{-1}^1u\p(\lambda u)\mathop{du}}{\int_{-1}^1\p(\lambda u)\mathop{du}} \\
& = \frac{1}{p\lambda}\left(\frac{(d+1)((p\lambda+c)^{d+2}-(c-p\lambda)^{d+2})}{(d+2)((p\lambda+c)^{d+1}-(c-p\lambda)^{d+1})}-c \right)  
\end{align*}
Note that, with our assumptions, \(\bary_{\p(\lambda\cdot)}([-1,1])>0\), so that we can write \(\beta_{\p}(\bbP^1,\ka)\) as: 
\[ \frac{p(d+2)((p\lambda+c)^{d+1}-(c-p\lambda)^{d+1})}{(\lambda p-c)(d
+2)((p\lambda+c)^{d+1}-(c-p\lambda)^{d+1})+(d+1)((p\lambda+c)^{d+2}-(c-p\lambda)^{d+2})} \]

\subsubsection{Fano simple principal \(\bbP^1\) bundles: comparison with Zhang and Zhou's work}
\label{sec_comp_ZZ}

It follows from classical results (see e.g. \cite[Lemma~2.1]{ZZ}) that, if \(L\) is an ample line bundle on \(B\) such that \(c_1(B)=rc_1(L)\) with \(r>1\), then the \(\bbP^1\)-bundle \(\bbP_{B}(\mathcal{O}\oplus L)\) is a Fano manifold. 
On the other hand, \(Y=\bbP_{B}(\mathcal{O}\oplus L)\) is a simple principal \(\bbP^1\)-bundle corresponding to the choices \(B=B_1= \SGr(n,2n+1)\), \(L_1=L\), \(p=1\). 
Note that, in this case, \(2\pi c_1(Y)\) is a compatible Kähler class, with \(\lambda=p=1\) and \(c=r\). 
However, we will show in Example~\ref{exa_not_comp_Fano} that such a Fano fibration may fail to be compatibly Fano. 

Assume now that \(Y\) is compatibly Fano. 
Let us compare our results with the results of Zhang and Zhou in \cite{ZZ}.
For \(\p(u)=(u+r)^d\) as above, we have 
\begin{align*} 
\beta_{\p}(\bbP^1,2\pi c_1(\bbP^1)) & = 
\frac{(d+2)((1+r)^{d+1}-(r-1)^{d+1})}{(1-r)(d+2)((1+r)^{d+1}-(r-1)^{d+1})+(d+1)((1+r)^{d+2}-(r-1)^{d+2})} \\
& = \left( 1-r + \frac{d+1}{d+2}\frac{(1+r)^{d+2}-(r-1)^{d+2}}{(1+r)^{d+1}-(r-1)^{d+1}} \right)^{-1}
\end{align*}
Up to the change in notations, this is precisely the quantity \(\beta_0\) in \cite{ZZ}, thus yielding a geometric interpretation of this quantity. 

Our expression for the greatest compatible Ricci lower bound is 
\begin{equation}
\beta^{\compatible}(Y,2\pi c_1(Y))=\min \left( \frac{r\beta(B,2\pi c_1(B))-1}{r-1}, \beta_0 \right)
\end{equation}
whereas the precise statement in \cite{ZZ} is 
\begin{equation}
\label{ZZ_formula}
\delta^{\alg}(Y) = \min \left(\frac{\delta^{\alg}(V)r\beta_0}{1+\beta_0(r-1)} , \beta_0\right) 
\end{equation}
where \(\delta^{\alg}(\cdot)\) denotes the algebraic delta invariant of a Fano variety. 
The inequality \(\beta^{\compatible}(Y,2\pi c_1(Y)) \leq \beta(Y,2\pi c_1(Y)) \leq \delta^{\alg}(Y) \) always hold. 
In fact, since \(Y\) admits a non-trivial \(\bbC^*\)-action, \(Y\) cannot be K-stable (although it could in principle be K-polystable) so \(\delta^{\alg}(Y)\leq 1\) and thus \(\beta(Y,2\pi c_1(Y))=\delta^{\alg}(Y)\). 
We will show in Example~\ref{exa_not_fiber} that the greatest compatible Ricci lower bound does not coincide with the greatest Ricci lower bound on a compatibly Fano example. 
On the other hand, it does coincide as soon as the minimum is computed by \(\beta_0\). More precisely, we have several equivalent properties:

\begin{prop}
The following are equivalent:
\begin{enumerate}
    \item \(\frac{r\beta(B,2\pi c_1(B))-1}{r-1}\geq \beta_0\)
    \item \(\frac{\delta^{\alg}(B)r\beta_0}{1+\beta_0(r-1)}\geq \beta_0\)
    \item \(\beta^{\compatible}(Y,2\pi c_1(Y)) = \beta(Y,2\pi c_1(Y))\)
\end{enumerate}
\end{prop}

\begin{proof}
First assume that (1) is satisfied. Then we have 
\[  \delta^{\alg}(B) \geq \beta(B,2\pi c_1(B)) \geq \frac{1+\beta_0(r-1)}{r} \]
reorganizing this inequality yields (2). 
By the formulas above we also have (3). 

Assume now that (2) is satisfied, then \(\beta(Y,2\pi c_1(Y)) = \delta^{alg}(Y)=\beta_0\). 
Assume by contradiction that \(\beta^{\compatible}(Y,2\pi c_1(Y))=\frac{r\beta(B,2\pi c_1(B))-1}{r-1} < \beta_0\). Then comparing with the above  inequalities, we must have \(\delta^{\alg}(B)> \beta(B,2\pi c_1(B))\). 
The latter is possible only if \(\beta(B,2\pi c_1(B))=1\). 
In this case, we have \(1=\frac{r\beta(B,2\pi c_1(B))-1}{r-1}<\beta_0\leq 1\), a contradiction. 
We have showed that (2) implies (3) (and (1)). 

Finally, assume (3), that is \(\beta^{\compatible}(Y,2\pi c_1(Y)) = \beta(Y,2\pi c_1(Y))\). If they are both equal to \(\beta_0\) then the formulas above imply (1) and (2). 
Assume now by contradiction that they are not equal to \(\beta_0\). 
We argue differently depending on the value of \(\delta^{\alg}(B)\). 
If \(\delta^{\alg}(B)\geq 1\), then \(\beta(B,2\pi c_1(B))\geq 1\), and as we have already seen, \(\beta^{\compatible}(Y,2\pi c_1(Y))=\beta_0\), a contradiction. 
If \(\delta^{\alg}(B)<1\), then \(\delta^{\alg}(B)=\beta(B,2\pi c_1(B))\). 
Now we consider the equality 
\[ \frac{r\beta(B,2\pi c_1(B))-1}{r-1} = \beta^{\compatible}(Y,2\pi c_1(Y)) = \beta(Y,2\pi c_1(Y)) = \frac{\delta^{\alg}(B)r\beta_0}{1+\beta_0(r-1)} \]
Using \(\delta^{\alg}(B)=\beta(B,2\pi c_1(B))\) and rearranging the terms, we deduce 
\[ \delta^{\alg}(B) = \frac{1+\beta_0(r-1)}{r} \]
and thus 
\[ \beta(Y,2\pi c_1(Y)) = \frac{\delta^{\alg}(B)r\beta_0}{1+\beta_0(r-1)} = \beta_0 \]
a contradiction again. 
\end{proof}

Our result shows that, in general, it is impossible to follow the smooth continuity method for Kähler-Einstein metrics while keeping compatibility of the metrics, and formula~\eqref{ZZ_formula} of \cite{ZZ} gives a hint of the interplay between the fiber and the basis needed to reach the end of the continuity method. 
It would be desirable to have a geometric interpretation of the behavior of metrics along the continuity method, or to find alternative techniques to compute the greatest Ricci lower bound of fibration by a differential geometric technique. 
An approach by conical twisted Kähler-Einstein metrics is sketched, using the Calabi ansatz in the case of Fano \(\bbP^1\)-bundles, in \cite[Appendix A]{ZZ}. 

\subsubsection{Examples}

\begin{exa}
Choose \(B_1=\bbP^{2}\), and \([\omega_1]=\frac{k}{3} 2\pi c_1(\bbP^{2})\) for some \(k\in \bbZ_{>0}\). 
Let \(p\in \bbZ_{>0}\) and consider the associated \(\bbP^1\)-bundle. 
For \(p=1\) and \(k=1\), we obtain the \(\bbP^1\)-bundle \(Y=\bbP_{\bbP^2}(\mathcal{O}\oplus \mathcal{O}(1))\), which is isomorphic to the blowup of \(\bbP^3\) at one point. 
By the above and Proposition~\ref{prop_Fano_toric_fiber} we have 
\[ \beta(Y,2\pi c_1(Y)) = \beta^{\compatible}(Y,2\pi c_1(Y)) = \beta_\p(\bbP^1,\ka) = \frac{14}{17} \]

For \(p=2\) and \(k=1\), or \(p=1\) and \(k=2\), we obtain the \(\bbP^1\)-bundle \(Y=\bbP_{\bbP^2}(\mathcal{O}\oplus \mathcal{O}(2))\), and 
\[ \beta^{\compatible}(Y,2\pi c_1(Y)) = \beta_\p(\bbP^1,\ka) = \frac{31}{43} \]
\end{exa}

\begin{exa}
\label{exa_not_fiber}
We now exhibit an example of compatibly Fano simple principal \(\bbP^1\)-bundle with  \(\beta^{\compatible}(Y, 2\pi c_1(Y)) \neq \beta_{\p}(\bbP^1,2\pi c_1(\bbP^1))\). 
Let us consider the previous example as a basis \(B_1=\bbP_{\bbP^2}(\mathcal{O}\oplus \mathcal{O}(1))\). This is a Fano threefold of index \(2\) which is not K-semistable. 
We consider the non-trivial Fano manifold \(Y\) defined as a simple principal \(\bbP^1\)-bundle by taking \(\lambda=1\), \(\ka = \frac{1}{2} 2\pi c_1(B_1)\), \(c=2\), \(p=1\). 
This fibration is indeed compatibly Fano since \(-1=\inf_{[-1,1]}p > -2\frac{14}{17}\). 
We compute, in this case, that 
\[ \beta_{\p}(\bbP^1,2\pi c_1(\bbP^1)) = \frac{50}{71} \]
and 
\[ \frac{2\beta(B_1,2\pi c_1(B_1))-1}{2-1} = \frac{11}{17} \]
thus 
\[ \beta^{\compatible}(Y,2\pi c_1(Y)) = \min\left\{\frac{50}{71},\frac{11}{17}\right\} = \frac{11}{17} \]

Using \cite{ZZ} we can furthermore compute 
\[ \beta(Y,2\pi c_1(Y)) = \min \left\{\frac{50}{71}, \frac{1400}{2057} \right\} = \frac{1400}{2057} \]
\end{exa}

\begin{exa}
\label{exa_not_comp_Fano}
Let us now exhibit an example of Fano simple principal \(\bbP^1\)-bundle which is \emph{not compatibly Fano}. 
Consider the \emph{odd symplectic Grassmannian} parametrizing \(n\)-dimensional isotropic subspaces in a \((2n+1)\)-dimensional complex vector space equipped with a skew-symmetric bilinear form of maximal rank. 
Since the dimension is odd, this is not a homogeneous space, but it has been studied a lot as an almost-homogeneous variety of Picard rank one \cite{Mihai_2007}, especially as a member of the family of horospherical Picard rank one varieties classified by Pasquier \cite{Pasquier_2009}. 
In particular, we know that \(\SGr(n,2n+1)\) is a \(\frac{n(n+3)}{2}\)-dimensional Fano manifold of Fano index \(n+2\) (see e.g. the table in \cite{Kanemitsu_2021}), which is not K-semistable \cite{Delcroix_KSSV}, and its greatest Ricci lower bound is 
\[ \beta(\SGr(n,2n+1), 2\pi c_1(\SGr(n,2n+1))) = \frac{2 ((2n+1)!)}{(n+2)(n! 2^n)^2} \sim \frac{4}{\sqrt{\pi n}}  \]
as proved by Hwang, Kim and Park in \cite{Hwang-Kim-Park_2023}. 
Recall that Fano index equal to \(n+2\) means that there exists a line bundle \(L\) on \(\SGr(n,2n+1)\) such that \((n+2)c_1(L) = c_1(\SGr(n,2n+1))\).
It follows from classical results (see e.g. \cite[Lemma~2.1]{ZZ}) that, if \(L\) is an ample line bundle on \(\SGr(n,2n+1)\) such that \(c_1(\SGr(n,2n+1))=rc_1(L)\) with \(r>1\), then the \(\bbP^1\)-bundle \(\bbP_{\SGr(n,2n+1)}(\mathcal{O}\oplus L)\) is a Fano manifold. 
On the other hand, \(Y=\bbP_{\SGr(n,2n+1)}(\mathcal{O}\oplus L)\) is a simple principal \(\bbP^1\)-bundle corresponding to the choices \(B=B_1= \SGr(n,2n+1)\), \(L_1=L\), \(p=1\). 
If \(Y\) were compatibly Fano we would have, for the compatible anticanonical class, \(c_1=r\), and 
\[-1=\inf_{[-1,1]}p > -r\beta(\SGr(n,2n+1), 2\pi c_1(\SGr(n,2n+1))) \]
Now for \(r\in \bbZ_{\geq 2}\), assume that \(n=r^3-2\) above, and consider the line bundle \(L=H^{r^2}\). 
Then we have 
\[ rc_1(L) = r^3 c_1(H) = (n+2) c_1(H) = c_1(B_1) \]
thus if \(Y\) is compatibly Fano we have \(c_1=r\) and  
\[-1 > -r\frac{2 ((2n+1)!)}{(n+2)(n! 2^n)^2} \sim -r\frac{4}{\sqrt{\pi(r^3-2)}} \sim -\frac{4}{\sqrt{\pi r}} \]
But the right hand side converges to \(0\) as \(r\) grows to \(+\infty\), a contradiction to the above inequality. 
As a consequence, for large \(r\), \(Y\) is Fano but not compatibly Fano. 
\end{exa}

\bibliographystyle{halpha}
\bibliography{logFibrations}
\end{document}